\documentclass[onefignum,onetabnum]{siamart171218}

\usepackage{amsfonts}
\usepackage{graphicx}
\usepackage{caption}
\usepackage{subcaption}
\usepackage{mathtools}
\usepackage{amssymb}
\usepackage{algorithmic}
\usepackage{array}
\usepackage{mathabx}
\usepackage[noadjust]{cite} 
\usepackage{comment}
\usepackage{amsopn}
\usepackage{placeins}

\ifpdf
\DeclareGraphicsExtensions{.eps,.pdf,.eps,.jpg}
\else
\DeclareGraphicsExtensions{.eps}
\fi
\newtheorem{example}{Example}[section]
\newtheorem{Remark}[theorem]{Remark}

\newcommand{\inner}[2]{\langle #1, #2 \rangle_F}

\headers{Circulant decomposition and Toeplitz type matrices}{Hariprasad M. and Murugesan Venkatapathi}

\title{Circulant decomposition of a matrix and the eigenvalues of Toeplitz type matrices}

\author{Hariprasad M.\and Murugesan Venkatapathi\thanks{Department of Computational \& Data Sciences, Indian Institute of Science, Bangalore.
		(\email{mhariprasadkansur@gmail.com}}, 
	\email{murugesh@iisc.ac.in}).}

\makeatletter
\newcommand*{\addFileDependency}[1]{
	\typeout{(#1)}
	\@addtofilelist{#1}
	\IfFileExists{#1}{}{\typeout{No file #1.}}
}
\makeatother

\DeclarePairedDelimiter\ceil{\lceil}{\rceil}

\begin{document}

	\DeclarePairedDelimiter\norm{\lVert}{\rVert}
	\maketitle
	
\begin{abstract}
We begin by showing that any $n \times n$ matrix can be decomposed into a sum of $n$ circulant matrices with periodic relaxations on the unit circle. This decomposition is orthogonal with respect to a Frobenius inner product, allowing recursive iterations for these circulant components. It is also shown that the dominance of a few circulant components in the matrix allows sparse similarity transformations using Fast-Fourier-transform (FFT) operations. This enables the evaluation of all eigenvalues of \emph{dense} Toeplitz, block-Toeplitz, and other periodic or quasi-periodic matrices, to a reasonable approximation in $\mathcal{O}(n^2)$ arithmetic operations. The utility of the approximate similarity transformation in preconditioning linear solvers is also demonstrated.
\end{abstract}
	
	
\section{Introduction}
Matrix decompositions, low-rank approximations, and projections on low complexity classes are some of the approaches useful in reducing the computation incurred with dense matrices \cite{Ambikasaran2013Dense}. Other approaches are also available in the case of sparse and effectively-sparse matrices, for example, in efficiently and accurately evaluating eigenvalues \cite{saad2011numerical, koslowski1993linear, paige1971computation, bunch2014sparse, ekstrom2018sparse, hariprasad2021chain}. Dense matrices that occur frequently in the numerical solution of eigenvalue problems such as Toeplitz, block-Toeplitz, and other matrices with periodicity in the diagonals are a class where further reduced computing may be possible. An efficient decomposition of the given matrix into circulant components i.e. circulant matrices with periodic relaxations on the unit circle, is shown to allow drastic reduction in such computations at a relatively small cost in the accuracy.

$\mathcal{O}(n^2)$ algorithms to evaluate the characteristic polynomial and its derivative were proposed a few decades ago for Toeplitz \cite{trench1989numerical}, block-Toeplitz \cite{trench1985eigenvalue, bini1988efficient} and Hankel matrices \cite{luk2000fast}. These algorithms were amalgamated with the Newton methods to evaluate only the eigenvalues of interest. On the other hand, the asymptotic behavior of spectra of Toeplitz and block-Toeplitz operators in the limit of large dimensions are well studied \cite{Gray1972Toeplitz, Boo1998Toeplitz, bottcher2005spectral, bottcher2012introduction,zhu2017asymptotic} and they can also be approximated by appropriate circulant matrices. Matrix less methods have been proposed for evaluating eigenvalues of certain classes of Toeplitz matrices
using asymptotic expansions \cite{bogoya2022fast}.

	
An approximation of a finite-dimensional Toeplitz matrix using a circulant matrix for speed up of operations was suggested decades ago \cite{chan1987asymptotic,arbenz1991computing}. Similarly, circulant preconditioners were constructed for Toeplitz matrices by minimising the Frobenius norm of the residue \cite{chan1988optimal}. Spectral preserving properties of this preconditioner were shown, and block circulant versions were also proposed and analysed \cite{jin2008survey}. They are expected to speed up iterative methods in regularization and optimization even though they can fail in special cases of poor approximation, and they may not be sufficiently accurate as a similarity transformation of the given Toeplitz matrix. This preconditioner based on minimizing the Frobenius norm of residue, is shown to be the first term in the circulant decomposition of the matrix discussed here. The relation between this single-term circulant approximation of a Toeplitz matrix, and its $symbol$, is highlighted in the appendix. The approach of minimizing the norm of the residue for a Toeplitz matrix, was extended to a generalized circulant preconditioner where only the absolute values of the entries preserve the periodicity \cite{Reichel2012circulant, Noschese2014circulant}. The full circulant decomposition of a matrix enhances and broadens the scope of such preconditioners and approximate similarity transformations, to matrices with any periodicity along the diagonals.

We first recall in \cref{remark-cycles} that any matrix can be decomposed into $n$ cycles that generate its $2n$-1 diagonals. We later show in \cref{theorem-circulant-similarity} and \cref{circD} that the decomposition of a matrix into $n$ circulant matrices with periodic relaxations on the unit circle, is equivalent to a decomposition of a similar matrix into such cycles. By including only the dominant cycles of the similar matrix, we include the dominant circulant components of the given matrix. The sparse similar matrix can be operated by a non-symmetric Lanczos algorithm for (block) tridiagonalization, and one can evaluate all eigenvalues of such dense matrices in $\mathcal{O}(n^2)$ arithmetic operations. Other relevant approaches for eigenvalues of sparse matrices have been reported as well \cite{saad2011numerical}.
	
Let $I_n$ be the identity matrix of dimension $n$, and $C$ be the permutation matrix corresponding to a full cycle.
	$$C = 
	\begin{bmatrix}
	0 & 1 \\
	I_{n-1} & 0 
	\end{bmatrix}_{ n \times n}.
	$$ 
	
	\begin{Remark} \label{remark-cycles}
		Any matrix A can be decomposed into $n$ cycles given by a power series in $C$ such that $A = \sum \limits_{k=0}^{n-1} \Lambda_kC^k$, where the Hadamard product $A \circ C^k = \Lambda_kC^k$, and $\Lambda_k$ are diagonal matrices. Entries supported on $C^k$ i.e. diagonal entries of $\Lambda_k$ in the above decomposition, are referred as the $k^{\text{th}}$ cycle of the matrix $A$.
	\end{Remark}

\begin{example}
The decomposition of a matrix into cycles, $A = \sum \limits_{j=0}^{n-1} \Lambda_jC^j$ for an order 3 magic square.
$$  \begin{bmatrix}
	\textcolor{green}{8} & \textcolor{blue}{1} & \textcolor{red}{6}\\
	\textcolor{red}{3} & \textcolor{green}{5} & \textcolor{blue}{7}\\
	\textcolor{blue}{4} & \textcolor{red}{9} & \textcolor{green}{2}
	\end{bmatrix}
	=
	\begin{bmatrix}
	\textcolor{green}8 & 0 & 0 \\
	0 & \textcolor{green}5 & 0\\
	0 & 0 & \textcolor{green}2
	\end{bmatrix} C^0 +
	\begin{bmatrix}
	 \textcolor{red}3 & 0 & 0 \\
	0 &  \textcolor{red}9 & 0\\
	0 & 0 &  \textcolor{red}6
	\end{bmatrix} C^1 +
	\begin{bmatrix}
	 \textcolor{blue}1 & 0 & 0 \\
	0 &  \textcolor{blue}7 & 0\\
	0 & 0 &  \textcolor{blue}4
	\end{bmatrix} C^2 $$ 
\end{example}
	
Let the permutation matrix given by a flipped identity matrix be
	$$
	J = 
	\begin{bmatrix}
	0  & \cdots&0 &0 & 1 \\
	0  & \cdots&0 &1 & 0 \\
	0  & \cdots&1 &0 & 0 \\
	\vdots  & .^{.^{.}}&0 &0 & 0 \\
	1 &  \cdots&0 &0 & 0 \\
	\end{bmatrix}_{ n \times n}.$$

For a sequence $\{x(p)\}_{0}^{n-1}$, its Discrete-Fourier-Transform (DFT) is given by $\{X(k)\}_{0}^{n-1}$, with $X(k) = \sum \limits_{p=0}^{n-1}x(p) e^{\frac{i 2 \pi pk}{n} }$. Eigenvalues of a circulant matrix $R$ is given by the DFT of the first row $\{R(0,j) \}_{j=0}^{n-1}$. The corresponding eigenvectors are given by the columns of the matrix $W$ with $W(p,q) = \frac{1}{\sqrt{n}}e^{-i \frac{2\pi pq}{n}}$, with $p,q=0,1,2\cdots n-1$.

\begin{Remark} \label{remark-circulant-eigen-properties}
$W^2 = CJ$, and any circulant matrix has an eigen decomposition $R = W \Lambda W^{\dagger}$. $R$ is also given by $R = W^{\dagger} \tilde{\Lambda} W$ where $\tilde{\Lambda} = CJ \Lambda J^T C^T$. Thus for $1 \leq k \leq n-1$, we have $\tilde{\Lambda}(k,k) = \Lambda(n-k,n-k)$ and $\Lambda(0,0) = \tilde{\Lambda}(0,0)$.
\end{Remark}
	
\begin{lemma} \label{theorem-circulant-similarity}
Given diagonal matrices $D_k$ with $D_k(q,q) = e^{i\frac{2 \pi k q}{n}}$ (for $ 0 \leq q  \leq n-1$) and any circulant matrix $R$ with eigenvalues given by a diagonal matrix $\tilde{\Lambda}$, the matrices $RD_k$ and $\tilde{\Lambda} C^k$ are similar.
\end{lemma}
\begin{proof}
Using $W$ for a linear transformation of $RD_k$,
\begin{align*}
WRD_kW^{\dagger} &= W W^{\dagger} \tilde{\Lambda} W D_k W^{\dagger} \\
& = \tilde{\Lambda} W D_k W^{\dagger} \\
&= \tilde{\Lambda} C^k. 
\end{align*}
The substitution $W D_k W^{\dagger} = C^k$, can be deduced by evaluating its $(p,l)^{\text{th}}$ entry.
\begin{align*}
(p,l)^{\text{th}}\text{ entry of } WD_kW^{\dagger} &= \frac{1}{n}\sum \limits_{q=0}^{n-1} e^{-i\frac{2 \pi pq}{n}} e^{i\frac{2 \pi kq}{n}} e^{i\frac{2 \pi ql}{n}} \\
&= \frac{1}{n}\sum \limits_{q=0}^{n-1} e^{i\frac{2 \pi (-p+k+l) q}{n}} \\
&= 1 \text{ when } p \equiv l+ k \mod n,\text{ and $0$ otherwise}.   
\end{align*}
Hence, the eigenvalues of $RD_k$ are the eigenvalues of the matrix $\tilde{\Lambda} C^{k}$. 
\end{proof}
	
Note that the matrix $\tilde{\Lambda} C^{k}$ is sparse and represents a single cycle, while the similar matrix $RD_k$ is dense. \cref{remark-cycles} recalls that any matrix has a decomposition into such cycles. The cycle decomposition of a matrix $A$ is equivalent to a circulant decomposition of a transformed similar matrix $WAW^{\dagger}$. Conversely, a circulant decomposition of the given matrix $A$ is equivalent to a cycle decomposition of the transformed similar matrix $WAW^{\dagger}$, as presented in the theorem below.
	
\begin{theorem} \label{circD}
Any $n \times n$ square matrix $A$ can be represented as a sum of $n$ circulant matrices with relaxations taking values from the $n^{\text{th}}$ root of unity. It is of the form
$A = \sum \limits_{k=0}^{n-1}R_kD_k$where $R_k$ is a circulant matrix and $D_k$ is a diagonal matrix with $D_k(q,q) = e^{i\frac{2 \pi kq}{n}}$ (with $ 0 \leq q  \leq n-1$).
\end{theorem}
	
\begin{proof}
Consider a matrix $B = WAW^{\dagger}$. Recalling its decomposition into cycles, and using the substitution $C^k = W D_k W^{\dagger}$ from the proof of \cref{theorem-circulant-similarity}, we have
$ B = \sum \limits_{k=0}^{n-1} \tilde{\Lambda}_k C^k = \sum \limits_{k=0}^{n-1} \tilde{\Lambda}_k W D_k W^{\dagger}$, where $\tilde{\Lambda}_k$ are diagonal matrices. 
Thus, the original matrix A is given by:
\begin{equation*}
A =  W^{\dagger}BW = \sum \limits_{k=0}^{n-1}   W^{\dagger} \tilde{\Lambda}_k W  D_k = \sum \limits_{k=0}^{n-1} R_k D_k .   
\end{equation*}
\end{proof}
	
\begin{example}
The decomposition of a matrix into circulant components, $A = \sum \limits_{k=0}^{n-1}R_kD_k$ for an order 3 magic square.
\begin{align*}
    \begin{bmatrix}
	8 & 1 & 6\\
	3 & 5 & 7\\
	4 & 9 & 2
	\end{bmatrix}
	&= \begin{bmatrix}
	5 & 4 & 6 \\
	6 & 5 & 4\\
	4 & 6 & 5
	\end{bmatrix} 
	\begin{bmatrix}
	1 & 0 & 0 \\
	0 & 1 & 0\\
	0 & 0 & 1
	\end{bmatrix} + 
	\begin{bmatrix}
	1.5-0.86i & 1.73i & -1.5-0.86i \\
	-1.5-0.86i & 1.5-0.86i & 1.73i\\
	1.73i & -1.5-0.86i & 1.5-0.86i
	\end{bmatrix} 
	\begin{bmatrix}
	1 & 0 & 0 \\
	0 & e^{i\frac{2\pi}{3}} & 0\\
	0 & 0 & e^{i\frac{4\pi}{3}}
	\end{bmatrix} \\ 
	& +
	\begin{bmatrix}
	1.5+0.86i & -1.73i & -1.5+0.86i \\
	-1.5+0.86i & 1.5+0.86i & -1.73i\\
	-1.73i & -1.5+0.86i & 1.5+0.86i
	\end{bmatrix} 
	\begin{bmatrix}
	1 & 0 & 0 \\
	0 & e^{i\frac{4\pi}{3}} & 0\\
	0 & 0 & e^{i\frac{2\pi}{3}}
	\end{bmatrix}.
\end{align*}
\end{example}

\textbf{Frobenius inner product:} Let $\inner{A}{B} = \sum \limits_{i,j} \widebar{A(i,j)} B(i,j)$ denote the Frobenius inner product of matrices $A$ and $B$. Note that any unitary transformation of the matrices preserves this inner product. It is also shown below that the circulant decomposition \ref{circD} is an orthogonal decomposition with respect to $\inner{.}{.}$. 

\begin{lemma}\label{finner}
    For any unitary matrices $U,V$ and matrices $A,B$,
    $\inner{A}{B} = \inner{UAV}{UBV}$.
\end{lemma}
\begin{proof}
    Let $\text{ }\widehat{}\text{ }$ denote a vector formed by concatenating the columns of a matrix such as $\widehat{A} = \begin{bmatrix}
    A(:,1) \\
    A(:,2) \\
    \vdots \\
    A(:,n)
    \end{bmatrix}$. Then $\inner{A}{B} = \widehat{A}^{\dagger} \hat{B}$. Let [ ] denote the expansion into a matrix of diagonal blocks, such as $[U] = \begin{bmatrix}
    U & 0 & \cdots & 0 \\
    0 & U & \cdots & 0\\
    0 & 0 & \ddots &  0 \\
    0 & 0 & \cdots & U 
    \end{bmatrix}_{n^2 \times n^2}$. Let $P$ be the permutation such that $P \widehat{A} = \widehat{A^T}$, and G=UAV, H=UBV. Then  
    \begin{align*}
    \widehat{G} = [V]P[U]\widehat{A} = T\widehat{A}.
    \end{align*}
     where $T = [V]P[U]$ is also a unitary matrix. Thus we have 
    \begin{align*}
        \inner{UAV}{UBV}=\inner{G}{H} &= \widehat{A}^\dagger T^\dagger T \widehat{B},\\
        &= \widehat{A}^\dagger \widehat{B}, \\
        &= \inner{A}{B}.
    \end{align*} 
\end{proof}

\begin{theorem}\label{fbasis}
In the decomposition, $A = \sum \limits_{k=0}^{n-1}R_kD_k$,  we have $\inner{D_i}{D_j} =n \delta_{i,j}$, (with $\delta_{i,j} = 1$ if $i=j$ and zero otherwise). Also, for $i \neq j$, $\inner{R_iD_i}{R_jD_j} = 0$. 
\end{theorem}
\begin{proof}
    It is easy to verify $\inner{D_i}{D_j} = n \delta_{i,j}$. From the transformation $A \to W A W^{\dagger}$, and by Lemma \ref{finner} we can show that,
    \begin{align*}
        \inner{R_iD_i}{R_jD_j} &= \inner{WR_iD_i W^\dagger}{WR_jD_jW^\dagger}, \\
        &=\inner{C^i \Lambda_i} {C^j \Lambda_j},\\
        &= 0 \text{ for } i \neq j.
    \end{align*}
\end{proof}

The \cref{fbasis} implies a Gram-Schmidt orthogonalization type of procedure to obtain the individual circulant matrices in the decomposition of \cref{circD}.

\begin{Remark}\label{circ_gram_schmidt}
\textbf{Recursive iterations for circulant components:} The following procedure gives the individual circulant components for a given matrix $A$.  
\begin{itemize}
    \item Initialize $A_0 = A$, 
    \item for $k$ = $0,1,2 \cdots n-2$
    \begin{itemize}
     \item $R_{k}(0,j) = \frac{1}{n} \left( \mathbf{1}^T( A_{k} \circ C^{j} )\mathbf{1} \right) \text{for } j = 0,1, \cdots n-1$, and $\mathbf{1}$ is a $n$-vector with entries all ones i.e. Average the entries of $A_{k}$ in the corresponding cycles to find the $n$ unknown entries of circulant $R_{k}$.
    \item $A_{k+1} = (A_{k} - R_{k})D_{-1}$.
    \end{itemize}
    \item $R_{n-1} = A_{n-1}$
\end{itemize}
\end{Remark}

The circulant matrix component $R_{k}$ of $A$ can also be evaluated using an inverse transformation of a cycle of $WAW^{\dagger}$, given by $W^{\dagger}(WAW^{\dagger} \circ C^{k})W$. Note that $R_0$ minimizes $\norm{A-R}_F$ for any circulant matrix $R$.  

\begin{lemma}\label{index_set}
	Given $\Lambda$ with diagonal entries $\{\lambda_j\}_{j=0}^{n-1}$,  let $R_1 = W \Lambda W^\dagger$ and $R_2 = W^\dagger \Lambda W$.  Then $|R_1(0,0)| = |R_2(0,0)|$ and $|R_1(0,j)| = |R_2(0,n-j)|$ for $1 \leq j \leq n-1$.  
\end{lemma}	

\begin{proof}
$R_1(0,j)$ is the $j$th coefficient of the inverse discrete Fourier transform of the diagonals of $\Lambda$. We have,
 \begin{align}
 |R_1(0,j)| &= |\sum \limits_{k=0}^{n-1} \lambda_k e^{i 2 \pi \frac{k j}{n}}| = |\sum \limits_{k=0}^{n-1} \lambda_k e^{-i 2 \pi \frac{k (-j)}{n}}| \\
 &= |\sum \limits_{k=0}^{n-1} \lambda_k e^{-i 2 \pi\frac{ k (n-j)}{n}}|.
 \end{align}   	
Note that the first row and first column of the matrices correspond to index zero of the discrete Fourier transform. So the claim follows.
\end{proof}	
We proceed further to derive the relationship between the frequencies in the varying entries along the diagonals of matrix, and the cycles of $WAW^{\dagger}$. The following definitions are relevant for this exercise.

\textbf{Weight of a cycle - $w_i$ :}
Let $B=\sum \limits_{i=0}^{n-1}R_iD_i$. The relative weight of circulant component $R_i$ in $B$ is $w_i = \frac{\norm{R_i}_F^2}{\sum \limits_i \norm{R_i}_F^2}$. When $B=WAW^{\dagger}$, note that $w_i$ also represents the relative weight of the cycle $i$ in $A$, and $\sum \limits_i w_i = 1$, $0 \leq w_i \leq 1$.

\textbf{Partial energy of a set of frequencies - $E_i$ :}
Let $A = \sum \limits_{i=0}^{n-1} \Lambda_i C^i$, and the discrete Fourier transform of the diagonal entries $\{\Lambda_i(j,j)\}_{j=0}^{n-1}$ be $\gamma^i_j$. Let $S_k = \{ a_1,a_2, \cdots a_k \}$ be a set of indices, with $\frac{\sum \limits_{j=1}^{k}(|\gamma^i_{a_j}|)^2}{\sum \limits_{j=1}^n (|\gamma^i_j|)^2 } = E_i$, where $0 \leq E_i \leq 1$.
We say that the frequencies given by $S_k$ have a partial energy $E_i$ on the cycle $i$.

Also note that $\delta^i_j$, the inverse discrete Fourier transform of the diagonal entries, have the corresponding index set $T_k = \{ b_1,b_2, \cdots b_k \}$ given by Lemma \ref{index_set} ($b_j = a_j$ if $a_j = 0$, else $b_j = n-a_j$), such that $\frac{\sum \limits_{j=1}^{k}(|\delta^i_{b_j}|)^2}{\sum \limits_{j=1}^n (|\delta^i_j|)^2 } = E_i$. 

\textbf{Relative magnitude of a set of cycles - $s$ :}
For a matrix $B = \sum \limits_{i=0}^{n-1} \Lambda_i C^i$, when the cycles corresponding to indices $J_k = \{j_1,j_2 , \cdots j_k \}$ have  $\frac{\sum \limits_{j \in J_k} \norm{\Lambda_j}_F^2}{\norm{B}_F^2} = s$, we say that the cycles $J_k$ have a relative magnitude $s$.

\begin{theorem}\label{th:dominance}
	When $E_i$ is the partial energy in a set of frequencies indexed by $S_k$ on the cycle $i$ of matrix $A$, the cycles of matrix $B=WAW^{\dagger}$ indexed by the corresponding $T_k$ have a relative magnitude $s = \sum w_i E_i$ where $\sum w_i = 1$ and $0 \leq w_i \leq 1$.
\end{theorem}
\begin{proof}
Consider the decomposition,
\begin{equation*}
	A = \sum \limits_{i=0}^{n-1} \Lambda_iC^i \text{ and }
	B = W (\sum \limits_{i=0}^{n-1} \Lambda_iC^i)  W^{\dagger}= \sum \limits_{i=0}^{n-1} R_iD_i.
\end{equation*}
By Parsevals theorem, $\norm{B}_F^2 = \sum \limits_i \norm{R_i}_F^2$, and in the decomposition $B = \sum \limits_{j = 0}^{n-1} \tilde{\Lambda}_j C^j$, cyclic diagonal entries $\tilde{\Lambda}_j$ are given by the discrete Fourier transform of the sequence  $\{ R_i(0,j)\}_{i=0}^{n-1}$ (and scaled by $\sqrt{n}$). Considering the ratio of entries corresponding to the cycles indexed by $T_k$,
\begin{align}
	\frac{\sum _{j \in T_k} \norm{\tilde{\Lambda}_j}_F^2}{\norm{B}_F^2} &= 
	n\frac{\sum _{j \in T_k} \sum \limits_{i=0}^{n-1}|R_i(0,j)|^2 }{\sum \limits_i \norm{R_i}_F^2},  \hspace{1cm} \text{using Parseval's theorem}\\
	&= 	n\frac{\sum \limits_{i=0}^{n-1}\sum _{j \in T_k}|R_i(0,j)|^2 }{\sum \limits_i \norm{R_i}_F^2}, \hspace{1cm} \text{changing the order of summation}
\end{align}

\begin{align}
	& = n\frac{\sum \limits_{i=0}^{n-1} \frac{\norm{R_i}_F^2}{n}   E_i}{\sum \limits_i \norm{R_i}_F^2}, \hspace{1cm} \text{by definition of $E_i$}\\
	& =  \sum \limits_{i=0}^{n-1} \left(\frac{\norm{R_i}_F^2}{\sum \limits_i \norm{R_i}_F^2}\right) E_i.
\end{align} 
	
Using $w_i = \frac{\norm{R_i}_F^2}{\sum \limits_i \norm{R_i}_F^2}$, we get $s = \sum \limits_{i=0}^{n-1} w_i E_i$.
\end{proof}

Note that $s \rightarrow 1$ as $E_i \rightarrow 1$. Thus, when a matrix $A$ has only $k \ll n$ dominant frequencies in the variation of entries along its diagonals, the similar matrix $B = WAW^{\dagger}$ is effectively sparse and has only $k$ dominant cycles, with the Frobenius norm of the other $n-k$ cycles being negligible.

\begin{corollary}
For a Toeplitz matrix 
$A = \begin{bmatrix}
a_0 & a_1 & a_2 &  \cdots & a_{n-1} \\
a_{-1} & a_0 & a_1 & a_2  & \ddots \\
a_{-2}& a_{-1} & a_0 & \ddots  & \ddots \\
\vdots & \ddots & \ddots    & \ddots & \ddots \\
a_{-(n-1)} & a_{-(n-2)}& \cdots & \cdots &a_0 \\
\end{bmatrix}$, by \cref{th:dominance}, the relative magnitude of the diagonal of $B = WAW^\dagger$ is $s^0=\sum \limits_i w_i E^0_i$, where $E^0_i$ includes only the zero-frequency or the average value of entries in the cycle $i$.
\begin{align}
    s^0 &= \sum \limits_i \frac{(n-i)|a_{-i}|^2 + i |a_{n-i}|^2}{\norm{A}^2_F} \frac{|(n-i) a_{-i} + i a_{n-i}|^2/n}{(n-i)|a_{-i}|^2 + i |a_{n-i}|^2} \\
    & =  \frac{\sum \limits_i |(n-i) a_{-i} + i a_{n-i}|^2}{n\norm{A}^2_F} \label{toeplitz_ratio}
\end{align}
In the above \eqref{toeplitz_ratio} note that for the cases $a_{-i} = a_{n-i}$, representing the circulant matrices, $s^0=1$ showing that only the diagonal of the corresponding $B$ has non-zero entries. Minimizing the first circulant component $s^0$ for a Toeplitz matrix using special pathological cases such as $a_{-i}=-ia_i/(n-i)$ where $E^0_i = 0$ for all $i > 0$, the lower bound on $s^0$ is 
\begin{equation*}
    s^0 \geq \frac{n|a_0|^2}{\norm{A}_F^2}
\end{equation*}
\end{corollary}

The above lower bound shows that the circulant matrix $R$ directly minimizing $\norm{A-R}_F$ need not be an effective approximation for a Toeplitz matrix in general, even in the limit of large $n$. Using expressions such as \eqref{toeplitz_ratio} for the other partial energies, one can show that when the entries of the Toeplitz matrix $\{a_{-(n-1)}, a_{-(n-2)},\cdots a_0, \cdots a_{(n-2)}, a_{(n-1)}\}$ are randomly chosen, partial energies $E^k_i$ and $E^{n-k}_i$ for frequencies $k \ll n$ are dominant corresponding to a small set of cycles in $B$. In other cases, depending on the set of given $2n-1$ entries, only the dominant cycles (representing the dominant partial energies) can be included in a sparse approximation of $B$ with a negligible residue.

\begin{corollary}
Given the set of entries $\{a_{-(n-1)}, a_{-(n-2)},\cdots a_0, \cdots a_{(n-2)}, a_{(n-1)}\}$ of a Toeplitz matrix, the distribution of the partial energies $E^k_i$ in the different frequencies indexed by $k=1,2,...n-1$, for a cycle $i$ is given by:
\begin{equation}
E^k_i = \frac{|a_{-i} - a_{n-i}|^2}{n ((n-i)|a_{-i}|^2 + i |a_{n-i}|^2 )} \left| \frac{\sin \frac{\pi (n-i+1)k }{n}}{ \sin \frac{\pi k}{n}} \right|^2
\end{equation}

\end{corollary}

\begin{corollary}
For a block Toeplitz matrix, we have constant matrices of a block size $m$ replacing the scalar entries $a_i$ and $a_{-i}$ in the Toeplitz matrix. Here, we have $m$ frequencies of interest given by the indices $S_m = \{n/m, 2n/m, 3n/m, \cdots n\}$. Correspondingly, we have the cycles of interest indexed by $T_m=\{n(m-1)/m, n(m-2)/m, \cdots 0\}$ in $B = WAW^\dagger$.


\end{corollary}

\section{Approximating eigenvalues using the circulant decomposition}

When a matrix $A$ has $k$ dominant frequencies in the variation of entries along its diagonals, the similar matrix $B = WAW^{\dagger}$ is effectively sparse and has $k$ dominant cycles, with the Frobenius norm of the other $n-k$ cycles being negligible (see \cref{th:dominance}). We can approximate $B$ by a sparse $\tilde{B} = \sum \limits_{i=1}^{k}\Lambda_{a_i} C^{a_i}$ by choosing the $k$ dominant cycles in $ a_i \in \{1,2 \cdots n \}$. This is followed by the application of an algorithm suited for the eigenvalues of a sparse matrix.
	
\subsection{Identifying the dominant circulant components}
\begin{figure}
		\parbox{7cm}{
	\includegraphics[width = 7 cm]{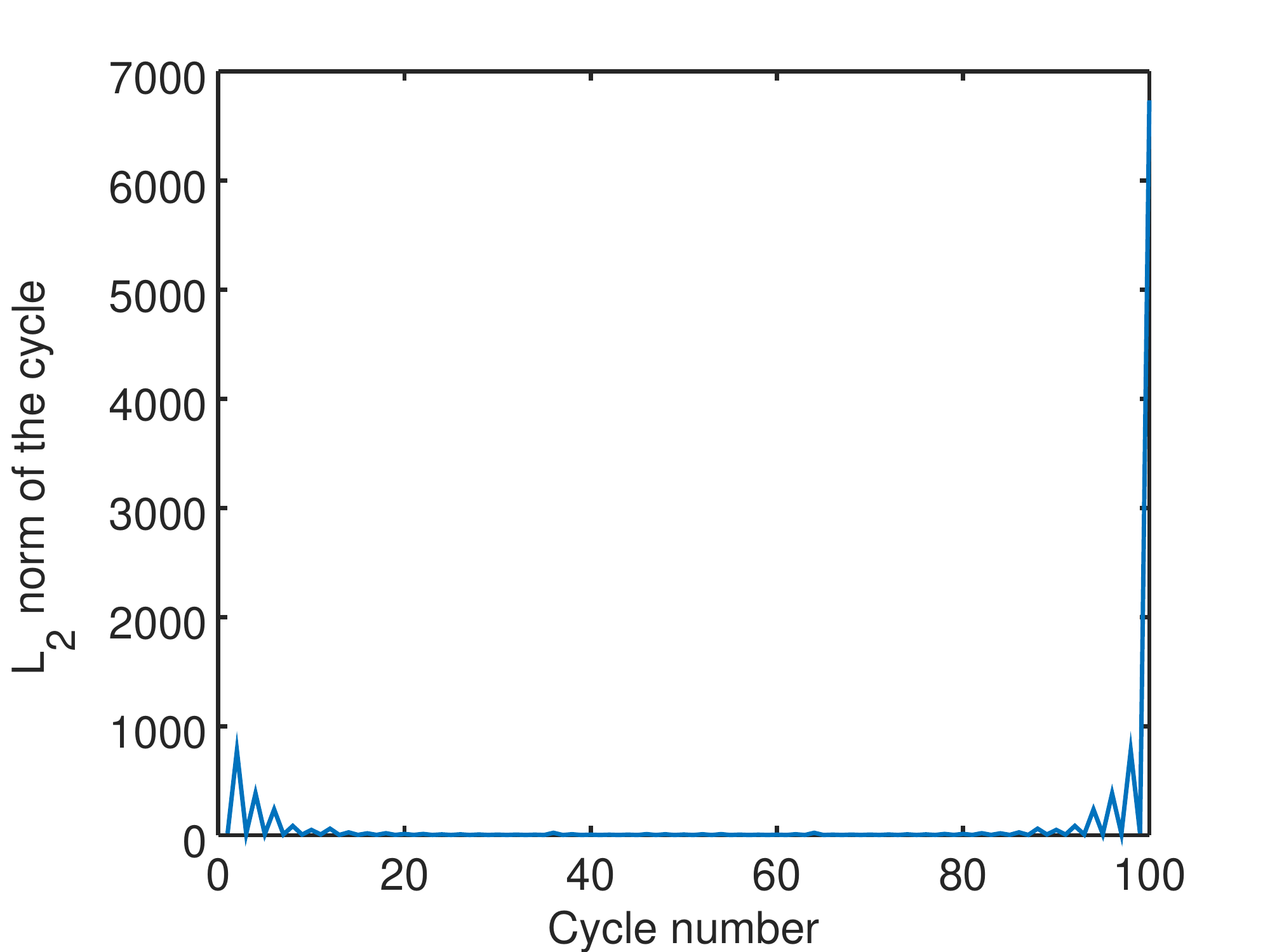}
		}
		\parbox{7cm}{
	\includegraphics[width = 7 cm]{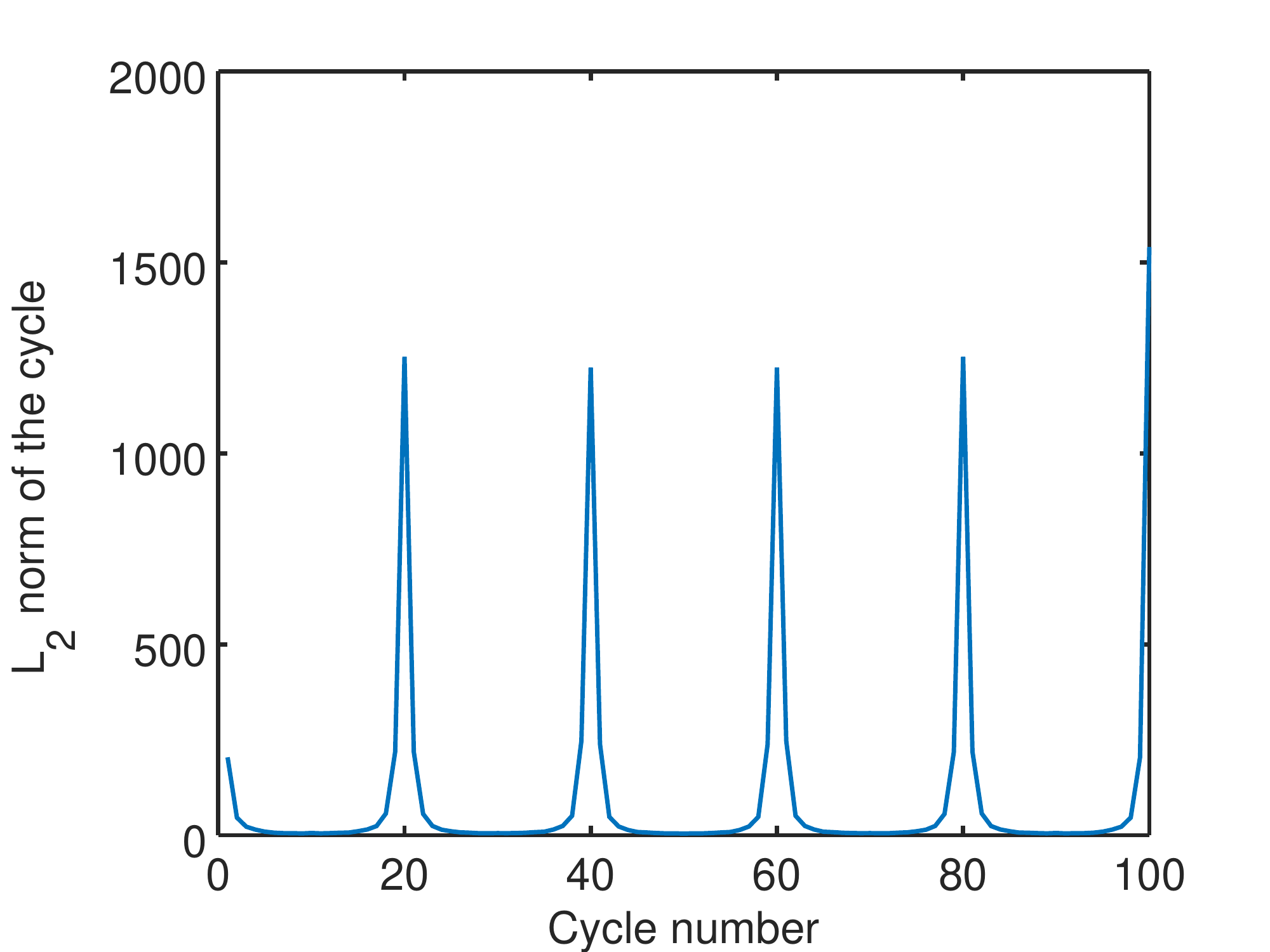}
		}
		\caption{$L_2$ norm of cycles of $WAW^{\dagger}$ numbered $0$ to $n-1$ where $n$=100. The dominant cycle numbering zero given by the diagonal is folded as number 100 in the plot. \textit{Left}: $A$ is a Toeplitz matrix with randomly chosen entries \textit{Right}: $A$ is a random block Toeplitz matrix of block-size 5.}\label{btc}
	\end{figure}

\begin{figure}
		\parbox{7cm}{
			\includegraphics[width = 7 cm]{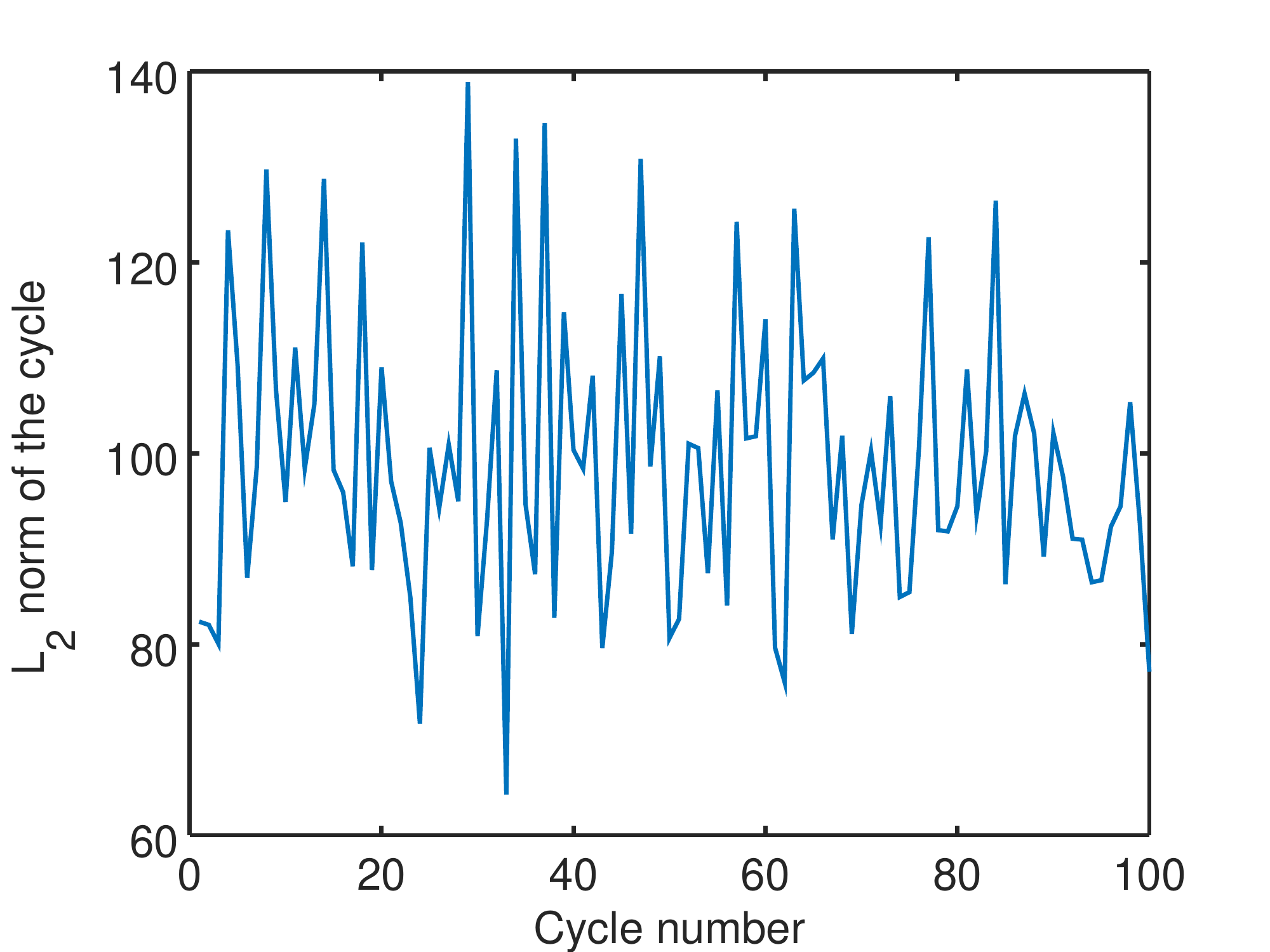}
		}
		\parbox{7cm}{
			\includegraphics[width = 7 cm]{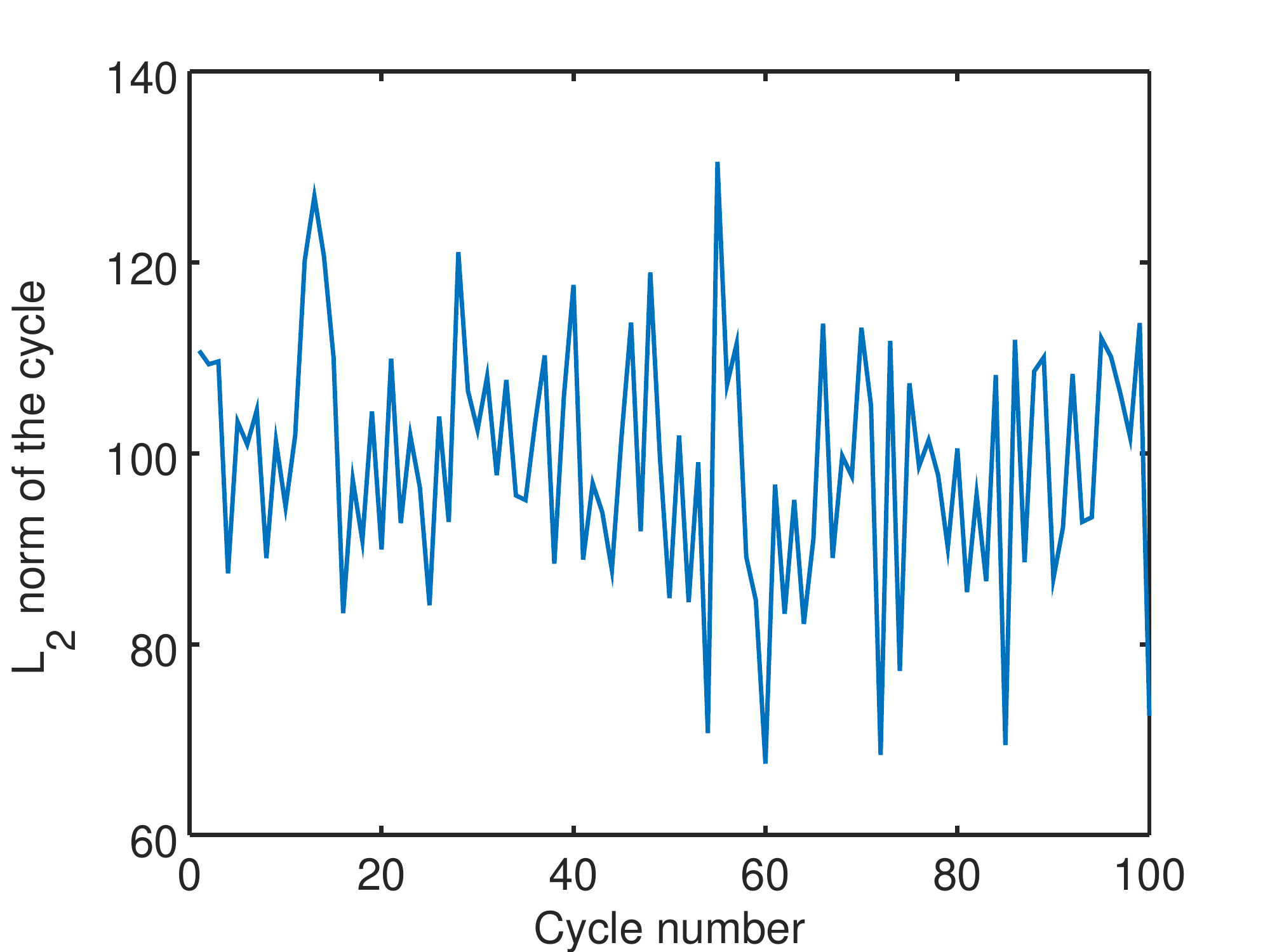}
		}
		\caption{The $L_2$ norm of cycles numbered $0$ to $n-1$ where $n$=100. The cycle numbered zero is folded as number 100 on the X-axis. \textit{Left}: Cycles of a random matrix A with entries from $\mathcal{N}(0,1)$. \textit{Right}: Cycles of $WAW^{\dagger}$ where entries of $A$ are given by $\mathcal{N}(0,1)$.}\label{wrc}
	\end{figure}

The similarity transformation by the matrix $W$ is shown to restrict the larger magnitudes to certain cycles
for the Toeplitz and block Toeplitz matrices (Figure \ref{btc}). However, the similarity transform fails to show such behavior for random matrices with no restriction on the entries (Figure \ref{wrc}), as expected.
	
Toeplitz matrices $A$, for example, have constant entries along the diagonals and the similarity transformation typically produces a very large magnitude for the cycle that represents a zero frequency (the main diagonal of the similar matrix $B$), and a few other cycles representing a low frequency of variations. On the other hand, a block Toeplitz matrix $A$ with a block size $m$ is $m$-periodic. Its similar matrix $B$ has dominant cycles given by integers $n/m$ and its multiples. Similarly, in the case of a quasi-periodic matrix where we have a random variation of $k$ frequencies along diagonals, the dominant cycles resemble a mixture of block-Toeplitz matrices of $k$ different block sizes. This is illustrated with examples in section \ref{sec-results}.

\subsection{Computing $WAW^{\dagger}$}

The $(p,q)$ entry of the matrix $B = W A W^{\dagger}$ is given by
\begin{align}
B(p,q) = \frac{1}{n}\sum \limits_{k=1}^{n} \sum \limits_{j=1}^{n} e^{-i\frac{2 \pi pk}{n}}A(k,j)e^{i\frac{2 \pi jq}{n}}.
\end{align}

\textit{Relation with fast Fourier transform} : If $\mathcal{F}(A)$ represents DFT of the columns of matrix $A$, we have 
\begin{align}
B = \frac{1}{n} \mathcal{F}(\mathcal{F}(A)^{\dagger})^{\dagger}.
\end{align}
If we denote the two dimensional DFT of the matrix $A$ by $\mathcal{F}^2(A)$, then we have $B = [\mathcal{F}^2(A)]/n$. Thus $WAW^{\dagger}$ is computed in $2n^2 \log n$ steps when the full linear transformation is required.

\textit{Evaluating only $k$ cycles of $WAW^{\dagger}$} : 
In the $n$-point fast-Fourier-transform (FFT) of a vector there are $\log_2n$ stages in its butterfly structure. When we require only $k$ frequency components of the given matrix at the final stage, the number of required points double every previous stage until $2^mk = n$ for some $m = \log_2\left( \frac{n}{k}\right)$. 
So the number of arithmetic operations $O_k$ for including $k$ frequency components is:
\begin{align*}
O_{k} &= \sum \limits_{j=0}^{m-1}2^jk + n\{\log_2n - \log_2 \left( \frac{n}{k} \right)\}, \\
O_{k} &= k (2^{m}-1) + n \log_2 k, \\
O_{k} &= (n-k) + n \log_2 k .
\end{align*}
The total number of required operations for $n$ vectors is $nO_k$, which is $\mathcal{O}(n^2)$ for a given $k$.
	
	\begin{algorithm}
		\textbf{Output :} {$\Lambda \longleftarrow$[ ]; a set of eigenvalues of given matrix.}\\
		\textbf{Initialization:} $W(p,q) \longleftarrow  \frac{1}{\sqrt{n}}e^{-i \frac{2\pi}{n} pq}$; Construct transformation matrix. \\
		\textbf{Sparsification :} $\tilde{B} \longleftarrow WAW^{\dagger} - \Delta$; Construct $\tilde{B}$ using the dominant cycles.\\
		\textbf{Reduced evaluation :} $\Lambda \longleftarrow \lambda_i\{\tilde{B}\}$; Eigenvalue algorithm for sparse matrices.
		\caption{Approximation of eigenvalues using a circulant decomposition}\label{approx}
	\end{algorithm}
	
\subsection{Error in approximations}
The similar matrix $B$ is reduced to a sparse matrix $\tilde{B}$ by selecting dominant cycles with the largest Frobenius norms. Let $\Delta = B-\tilde{B}$, and $\lambda_i$ be the eigenvalues of $B$, and $\tilde{\lambda}_i$ be the eigenvalues of $\tilde{B}$. Let $B= X \Lambda X^{-1}$. From Bauer-Fike theorem,
\begin{equation}
    |\lambda_i - \tilde{\lambda}_i| \leq \kappa(X) \norm{\Delta}_2
\end{equation}

Here $\kappa(X) = \norm{X}\norm{X^{-1}}$ is the condition number of the eigenvector matrix $X$. So the eigenvalues are better approximated when the eigenvector matrix $X$ is well conditioned, and $\norm{\Delta}_2$ is minimized. Similarly, the relative error can be bound when $B$ is non singular and diagonalizable.
\begin{align}
	\frac{ |\lambda_i - \tilde{\lambda}_i| }{ |\lambda_i| } & \leq \kappa(X)  \norm{B^{-1}\Delta}_2.
\end{align}

\subsection{Positive definiteness of $\tilde{B}$}

It can be shown that for a positive definite $A$, the diagonal matrix with the diagonal entries of $B$ is positive definite \cite{cai2005generalization}.
But the matrix $\tilde{B}$ with other cycles included, need not be positive definite. Here we provide a sufficient condition for the matrix $\tilde{B}$ to be positive definite.

\begin{theorem}
When the approximation $\tilde{B}$ includes the diagonal of $B$ along with $k$ symmetric cycles indexed by $a_i \in \mathcal{T}$, the approximated matrix is positive definite if $\frac{\sqrt{B(i,i)B(a_j,a_j)}}{k} \geq B(i,a_j)$ for all $i, a_j$. 
 \end{theorem}
\begin{proof}
By using the expansion of $u^\dagger \tilde{B} u$ with constraint $\norm{u} = 1$. The expansion can be divided into two parts
\begin{align}
    u^\dagger \tilde{B} u & = \sum \limits_i |u_i|^2 B(i,i) + \sum \limits_{p,q \in \mathcal{T}}  u_p \widebar{u_q} B(p,q) +  u_q \widebar{u_p} \widebar{B(p,q)}, \\
    u^\dagger \tilde{B} u & \geq \sum \limits_i |u_i|^2 B(i,i) - \sum \limits_{p,q \in \mathcal{T}} 2 |u_p| |u_q| |\text{Re}(B(p,q))|, \\
    u^\dagger \tilde{B} u & \geq \sum \limits_{p,q \in \mathcal{T}}\left(  \frac{(|u_p|^2 B(p,p) + |u_q|^2 B(q,q))}{|\mathcal{T}|} - 2 |u_p| |u_q| |\text{Re}(B(p,q))| \right).  
\end{align}

From the last equality, for every $p,q \in \mathcal{T}$ if
\begin{align}
    \frac{(|u_p|^2 B(p,p) + |u_q|^2 B(q,q))}{|\mathcal{T}|} - 2 |u_p| |u_q| |\text{Re}(B(p,q))| \geq 0, 
\end{align}
then it satisfies a sufficient condition for the matrix $\tilde{B}$ to be positive definite. With $\theta = \left|\frac{u_p}{u_q} \right|$, this implies
\begin{align}
    \frac{ \theta B(p,p) + \frac{1}{\theta} B(q,q)}{|\mathcal{T}|} \geq 2|\text{Re}(B(p,q))|
\end{align}
 
Now minimizing $f(\theta) = \theta B(p,p) + \frac{1}{\theta} B(q,q)$ w.r.t $\theta$, we find the minimum value to be $2\sqrt{B(p,p)B(q,q)}$, giving us the sufficient condition,
\begin{align}
    \frac{\sqrt{B(p,p)B(q,q)}}{|\mathcal{T}|} \geq |\text{Re}(B(p,q))|
\end{align}

On the other hand, $\frac{B(i,i) + B(j,j)}{2} < | \text{Re}(\tilde{B}(i,j))|$ for any $i,j$ implies a vector $u = ae_i + b e_j$ with suitable $a,b$ such that $u^\dagger B u < 0$. So we note that $\frac{B(i,i) + B(j,j)}{2} > | \text{Re}(\tilde{B}(i,j))|$ for all $i,j$ when $B$ is positive definite.   

\end{proof}

\section{Numerical Results} \label{sec-results}

\subsection{Eigenvalue evaluations}
We begin by qualitatively highlighting the advantages of the suggested sparsification over a direct sparsifier, for matrices with periodic properties. Later, we present quantitative results of the errors in the similarity transformation and eigenvalue approximations of Toeplitz, block Toeplitz and quasi-periodic matrices. Random matrices with $\mathcal{N}(0,1)$ as entries of the matrix or its blocks were used for corresponding numerical experiments on Toeplitz and block Toeplitz matrices. It should be noted that in applications where the matrices are not random, one can expect even better results. The results show that including a few dominant cycles provides us lower relative errors for eigenvalues of such matrices, compared to the single-term circulant approximation. Also, while the relative errors may be notable for the smaller eigenvalues, their absolute errors are very small, as expected.

\subsubsection{Comparison with direct sparsifiers}
\begin{figure}
		\parbox{7cm}{
			\includegraphics[width = 7 cm]{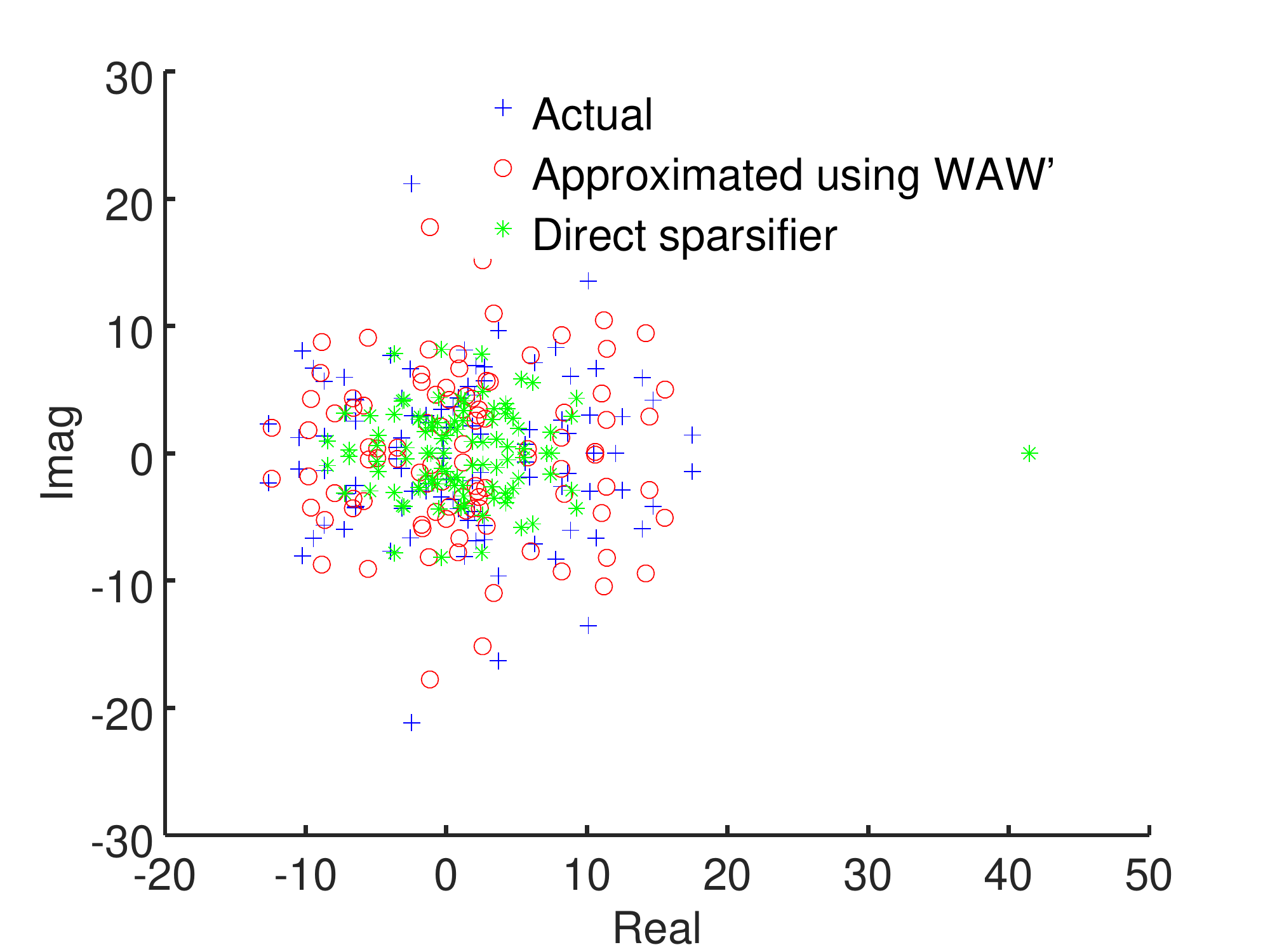}
		}
		\parbox{7cm}{
			\includegraphics[width = 7 cm]{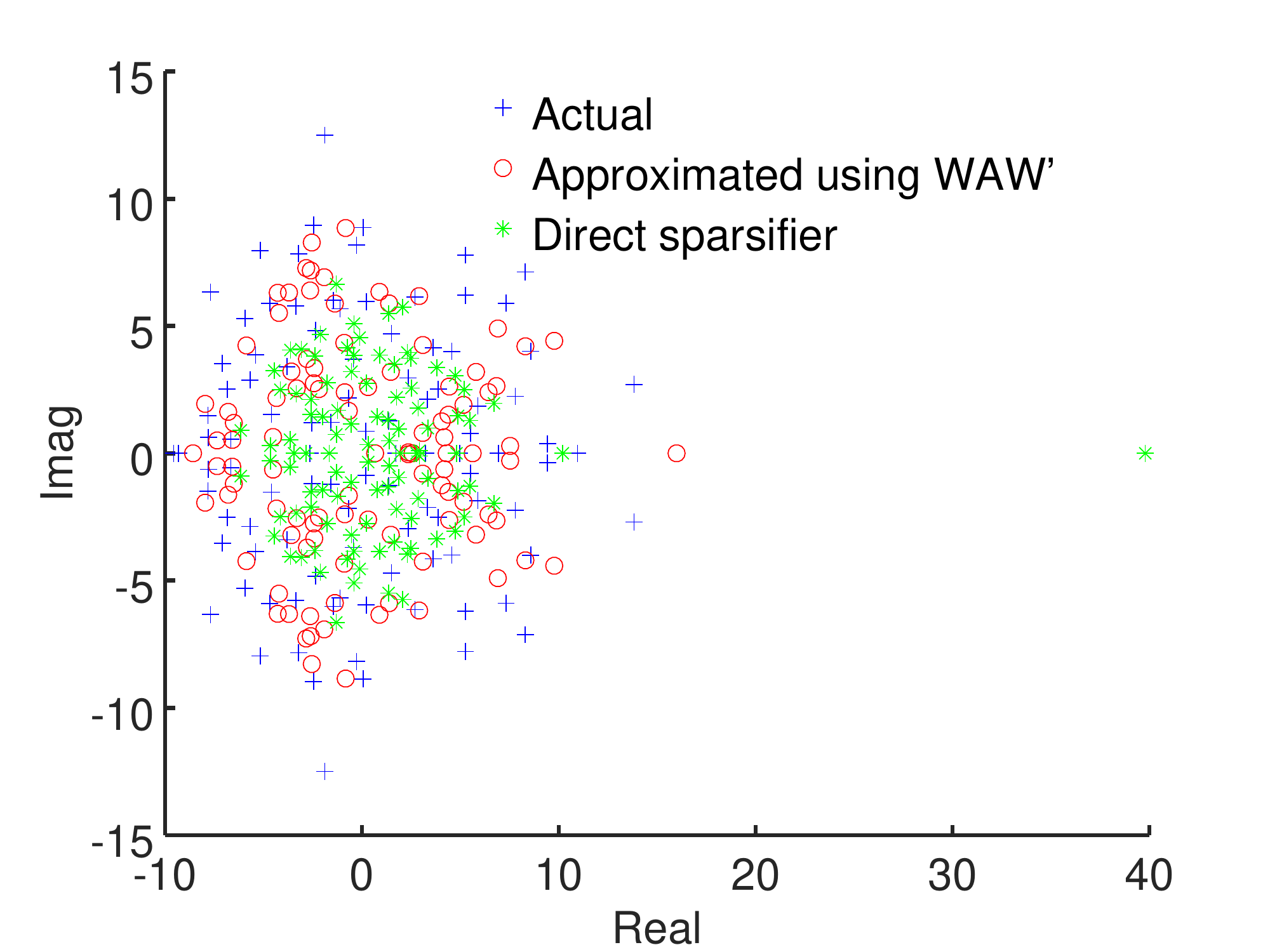}
		}
		\caption{A qualitative comparison of eigenvalue approximations using the direct \cite{achlioptas2007fast} and the circulant sparsifiers for identical number of non-zero entries in the resulting matrix. The smaller magnitudes of the former show that the information lost by the direct sparsifier is significantly larger compared to that of the suggested computationally efficient sparsification using the dominant circulant components. \textit{Left}: A random Toeplitz matrix. \textit{Right}: A random block Toeplitz matrix of block size 5.}\label{bte}
	\end{figure}

The eigenvalues of the approximated Toeplitz and block Toeplitz matrices using $WAW^{\dagger}$ and a direct sparsifier \cite{achlioptas2007fast}, along with the actual eigenvalues are shown in Figure \ref{bte}. In general, the eigenvalues with the direct sparsifier are smaller in magnitude i.e. more centered in the plot, compared to the eigenvalues with the proposed sparsifier in the frequency domain. Intuitively, the information lost by the direct sparsifier by thresholding some entries in the matrix to zero, is larger compared to the information lost in the frequency domain by removing the same number of entries as cycles of the matrix $WAW^{\dagger}$.
	
\subsubsection{Toeplitz and block-Toeplitz matrices}
\begin{figure}
		\parbox{7cm}{
			\includegraphics[width = 7 cm]{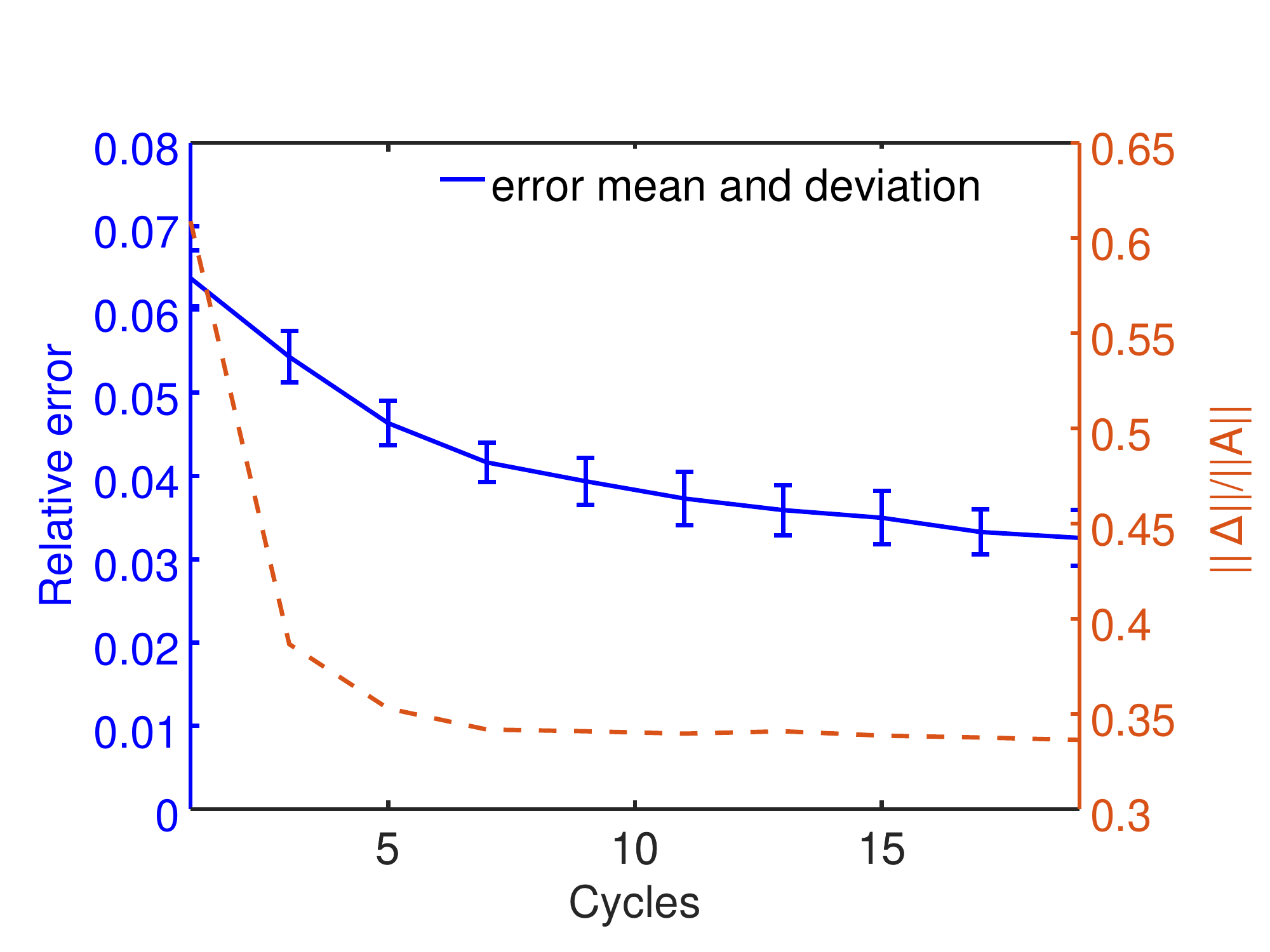}
		}
		\parbox{7cm}{
			\includegraphics[width = 7 cm]{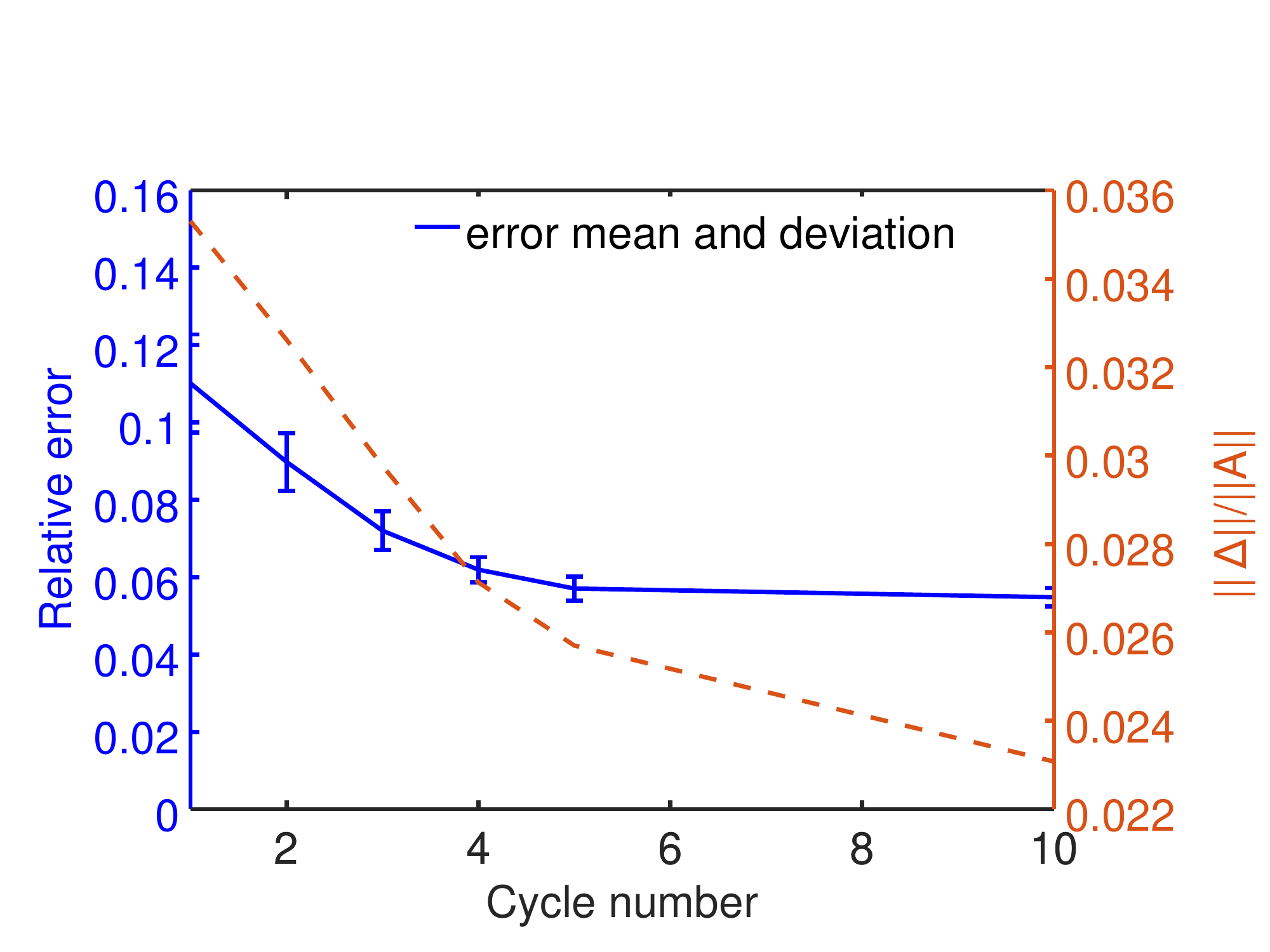}
		}
			\caption{Average and deviation of relative error in the approximated eigenvalues with increasing number of cycles considered. Corresponding relative errors in the similarity transformation are given by the dashed line and Y-axis on the right. \textit{Left}: Toeplitz matrix of dimension $n$=1000.  \textit{Right}: Block Toeplitz matrix of block size 5 and dimension $n$=1000.}\label{errors_toeplitz}
	\end{figure}

\begin{figure}
		\parbox{7cm}{
			\includegraphics[width = 7 cm]{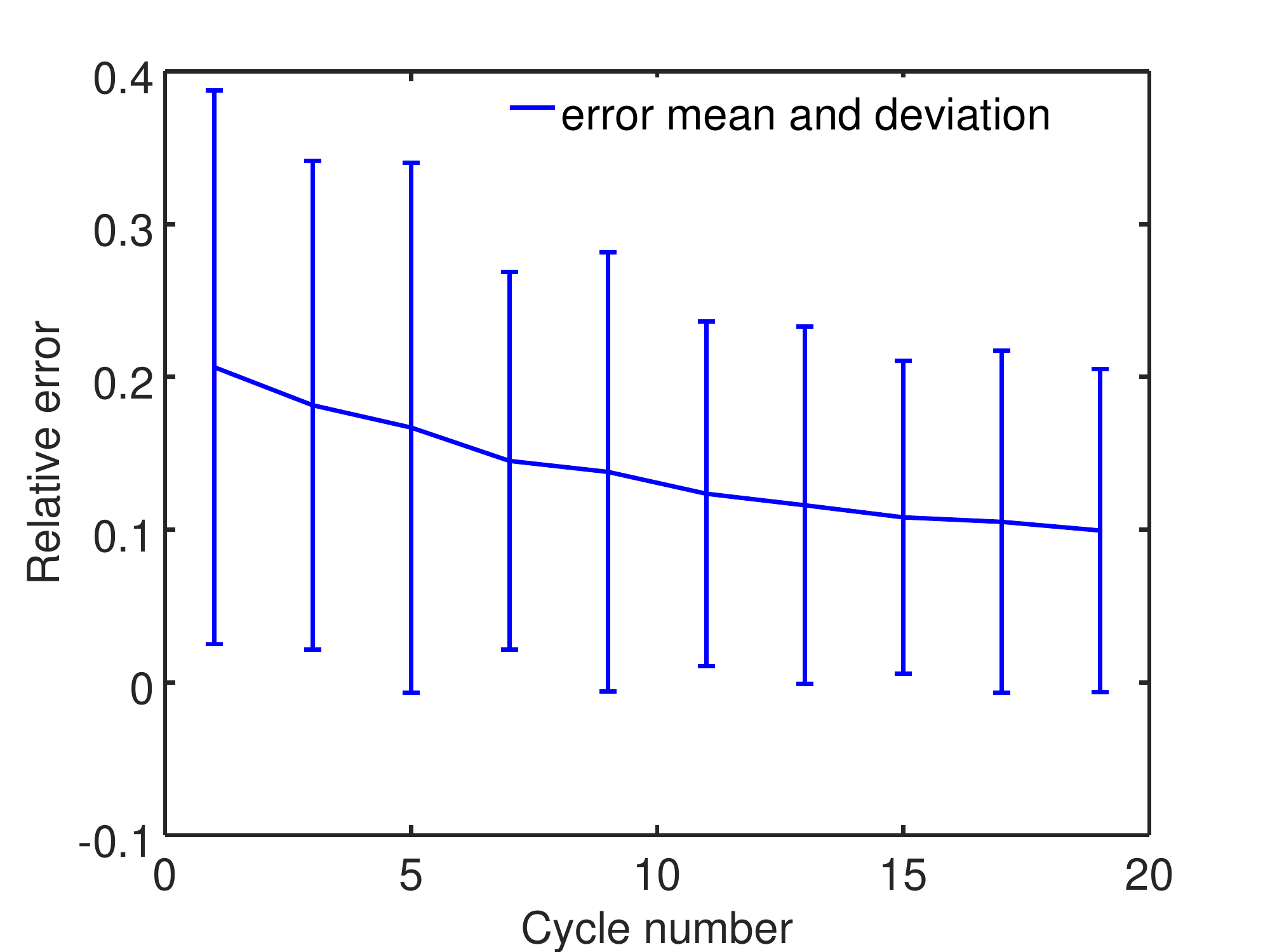}
		}
		\parbox{7cm}{
			\includegraphics[width = 7 cm]{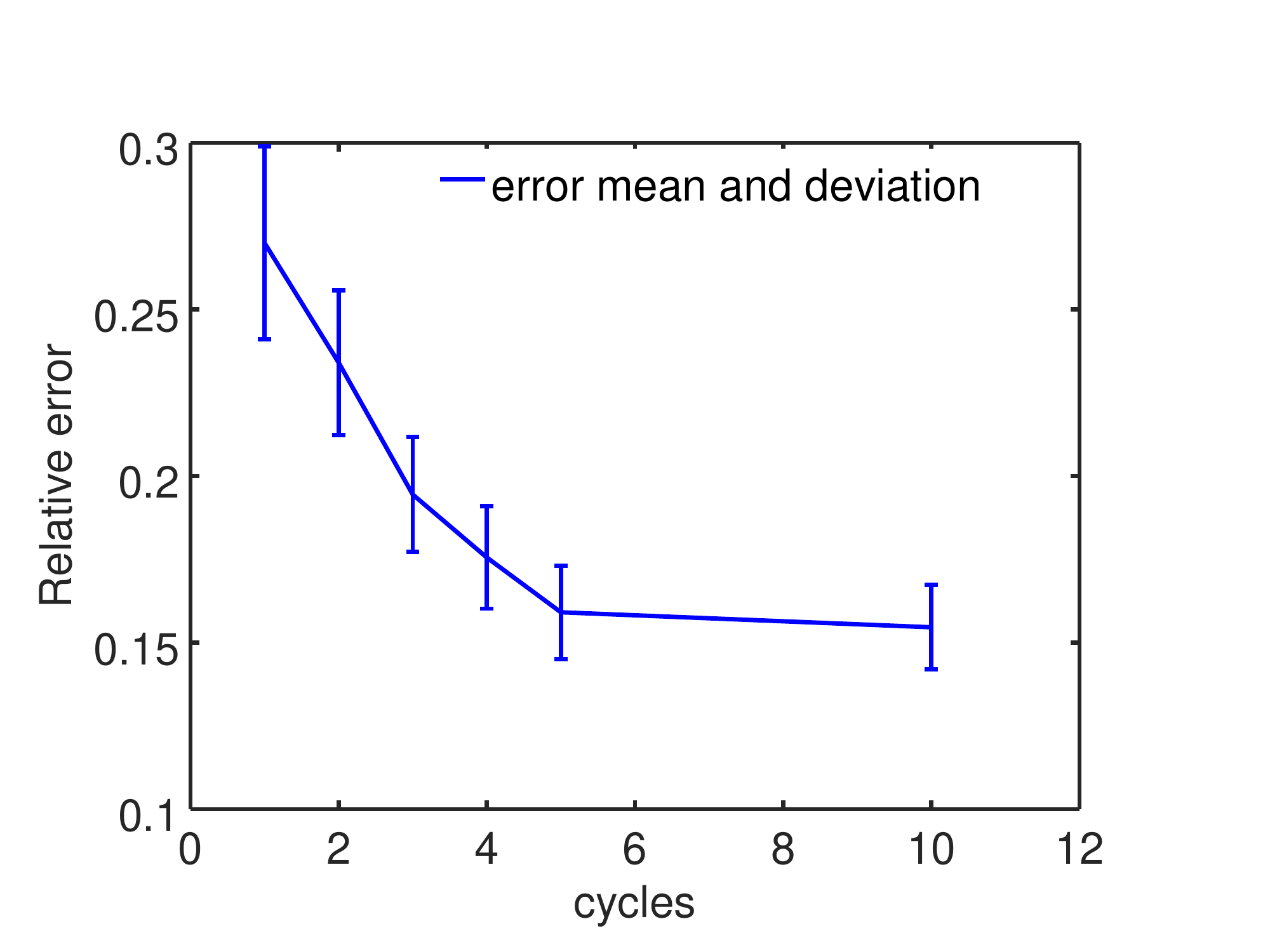}
    }
    \caption{Mean of the average relative error and its deviation in the approximated eigenvalues with the increasing number of included cycles; evaluated using 1000 random matrices of dimension $n$=100. Note that the errors further reduce significantly with the increasing dimensions $n$ as shown in Figure \ref{last_figure}. \textit{Left}: Toeplitz matrices. \textit{Right}: Block Toeplitz matrices of block size 5.}\label{average_errors_toeplitz}
	\end{figure}

The average relative error and standard deviation of the errors of all the $n$ eigenvalues, decrease with an increase in the number of cycles in the approximation of eigenvalues of a Toeplitz matrix, as shown in Figure \ref{errors_toeplitz}. The cumulative relative error of the approximate similarity transformation in terms of the Frobenius norms, $\norm{\Delta}/\norm{A}$, is also plotted in these figures. Also, the mean of the average relative error and its deviation over 1000 matrices in Figure \ref{average_errors_toeplitz} show reasonably small errors, and better results are expected for block Toeplitz and other periodic matrices.
	
Block-Toeplitz matrices with random entries from $\mathcal{N}(0,1)$ were used for the numerical experiments. The average relative error and standard deviation of the errors decrease with an increase in the number of cycles in the approximation as shown in the example in Figure \ref{errors_toeplitz}. The mean of the average relative error and its deviation over 1000 matrices are shown in Figure \ref{average_errors_toeplitz}. Note that inclusion of a few cycles may be sufficient to achieve a reasonable accuracy in approximating all the eigenvalues for these matrices. The improvements in the accuracy of the eigenvalue estimation for a given $k$, with the increase in the size of the matrix, is also illustrated in the Figure \ref{last_figure}.

\subsubsection{Other periodic and quasi-periodic matrices}
	
		\begin{figure}[h]
		\parbox{7cm}{
			\includegraphics[width = 7 cm]{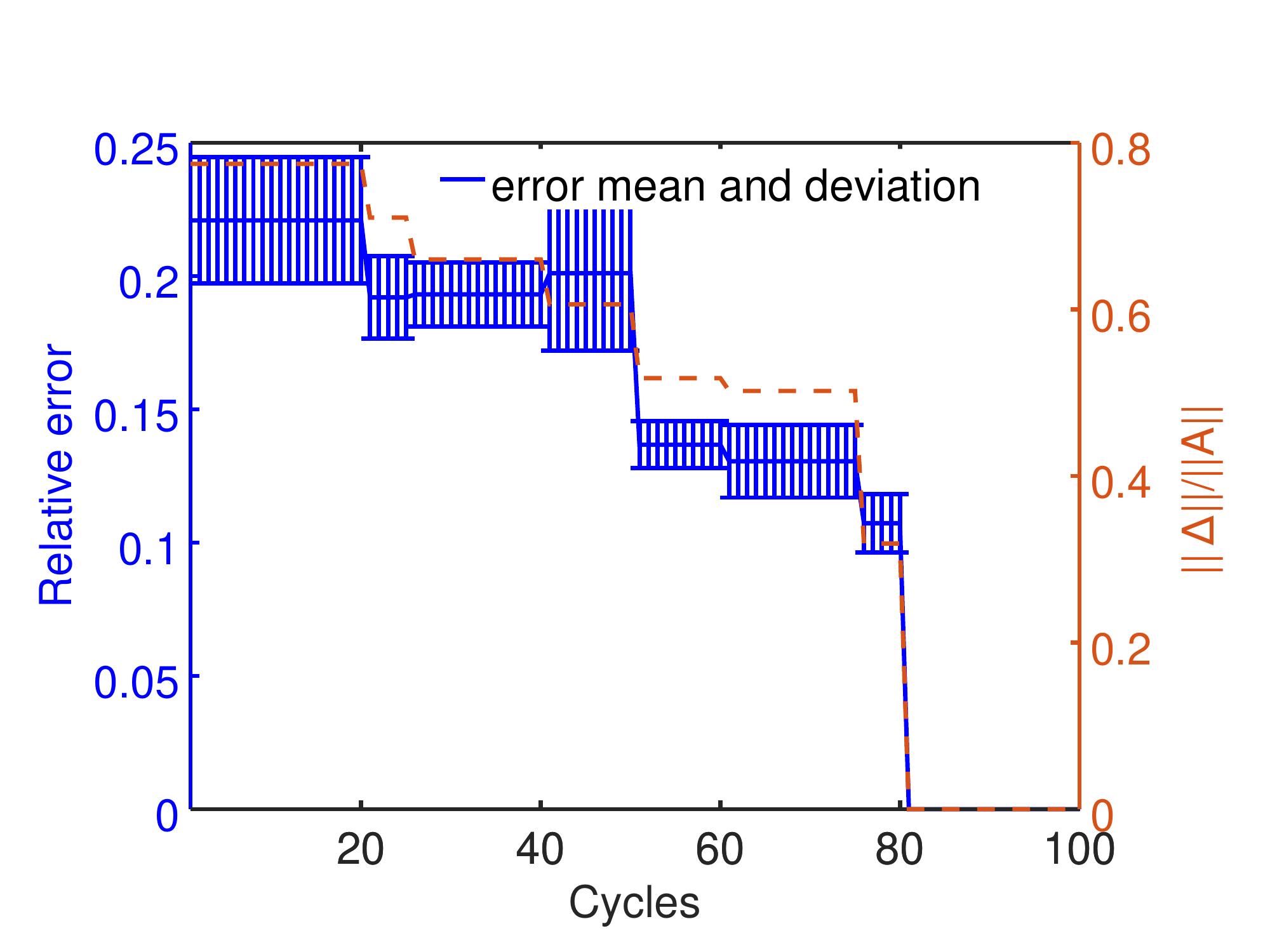}
		}
			\parbox{7cm}{
				\includegraphics[width = 7 cm]{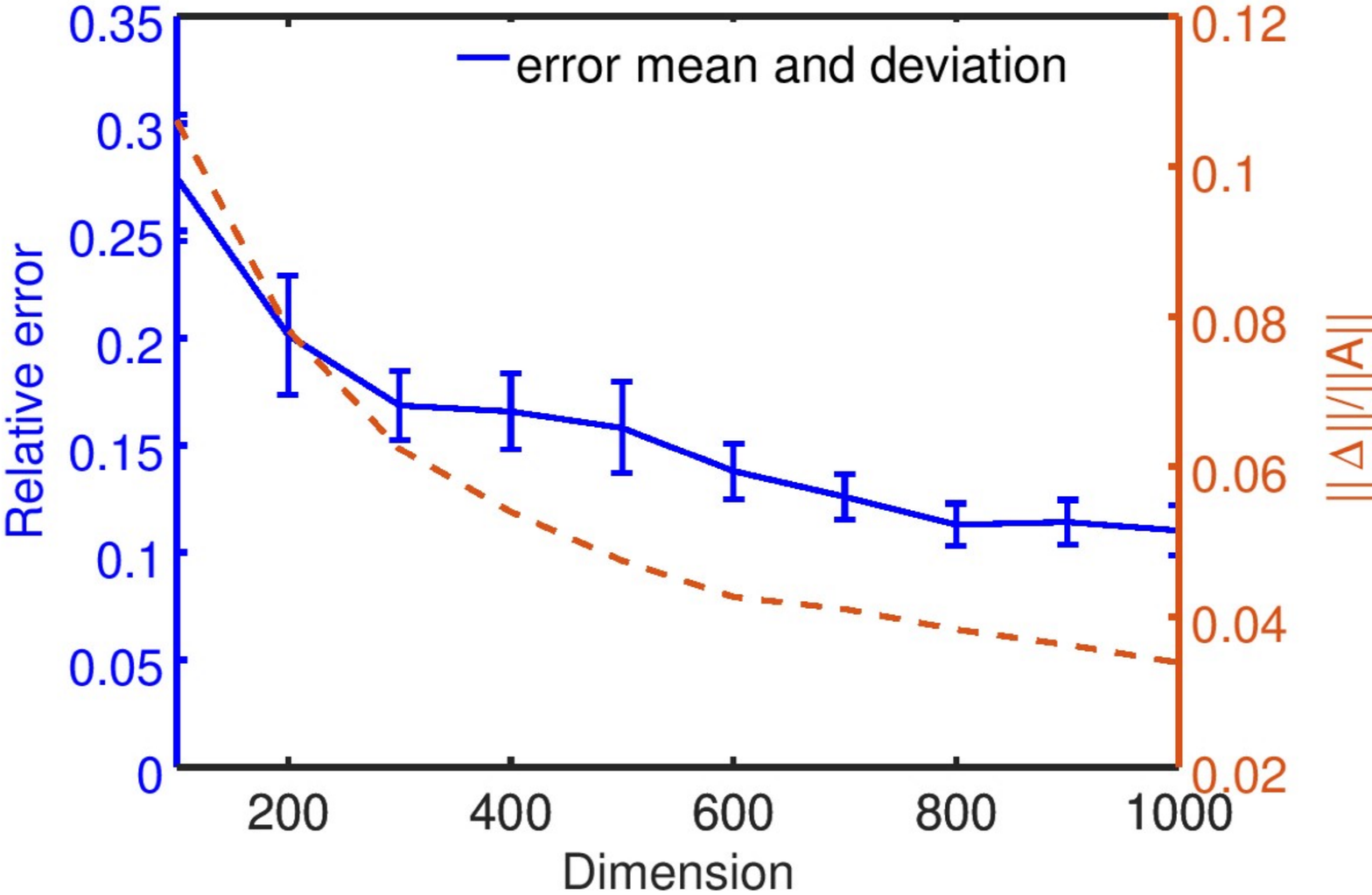}
		}
		\caption{\textit{Left}: Average and deviation of relative error in the approximated eigenvalues of a random periodic matrix of dimension $n$=100 with the included cycles. Periods 4, 5, 10 were chosen randomly along the diagonals with a uniform probability. In the suggested algorithm, cycles numbering around 20, 40, 50, 60 and 80 are preferentially included in the sparse similarity transformation. \textit{Right}: Average and deviation of relative error in the approximated eigenvalues using only two cycles ($C^0$, $C^{\lfloor \frac{n}{2} \rfloor}$), of block Toeplitz matrices of block-size five, for dimensions $n$ ranging from 100 to 1000. Corresponding relative errors in the similarity transformation are given by the dashed line and Y-axis on the right.}\label{last_figure}
	\end{figure}
	
The above properties of Toeplitz and block-Toeplitz matrices can describe the properties and the efficacy of sparsification of a quasi-periodic matrix as well. The dominant components for such matrices becomes the corresponding sizes of periodic blocks in the matrix; thus displaying the characteristics similar to that of a block-Toeplitz matrix. It can be used to include only the dominant cycles of the similar matrix $B$. Average relative error and deviation in the relative errors are plotted with the cycle number for a random periodic matrix in Figure \ref{last_figure}. Here, a uniform probability of periods $m$=4, 5 and 10 were used to generate the entries along the diagonals for a matrix of dimension $n$=100. These result in dominant cycles in the similar matrix $B$ that are multiples of the corresponding integers $n/m$ given by 25, 20 and 10 respectively, with their common multiples being even more significant. Note that the error in Figure \ref{last_figure} reduces mostly for corresponding cycles, as in the case of a mixture of block Toeplitz matrices.

\subsection{Preconditioning linear systems}
\begin{figure}[h]
	\parbox{7cm}{
		\includegraphics[width = 7 cm]{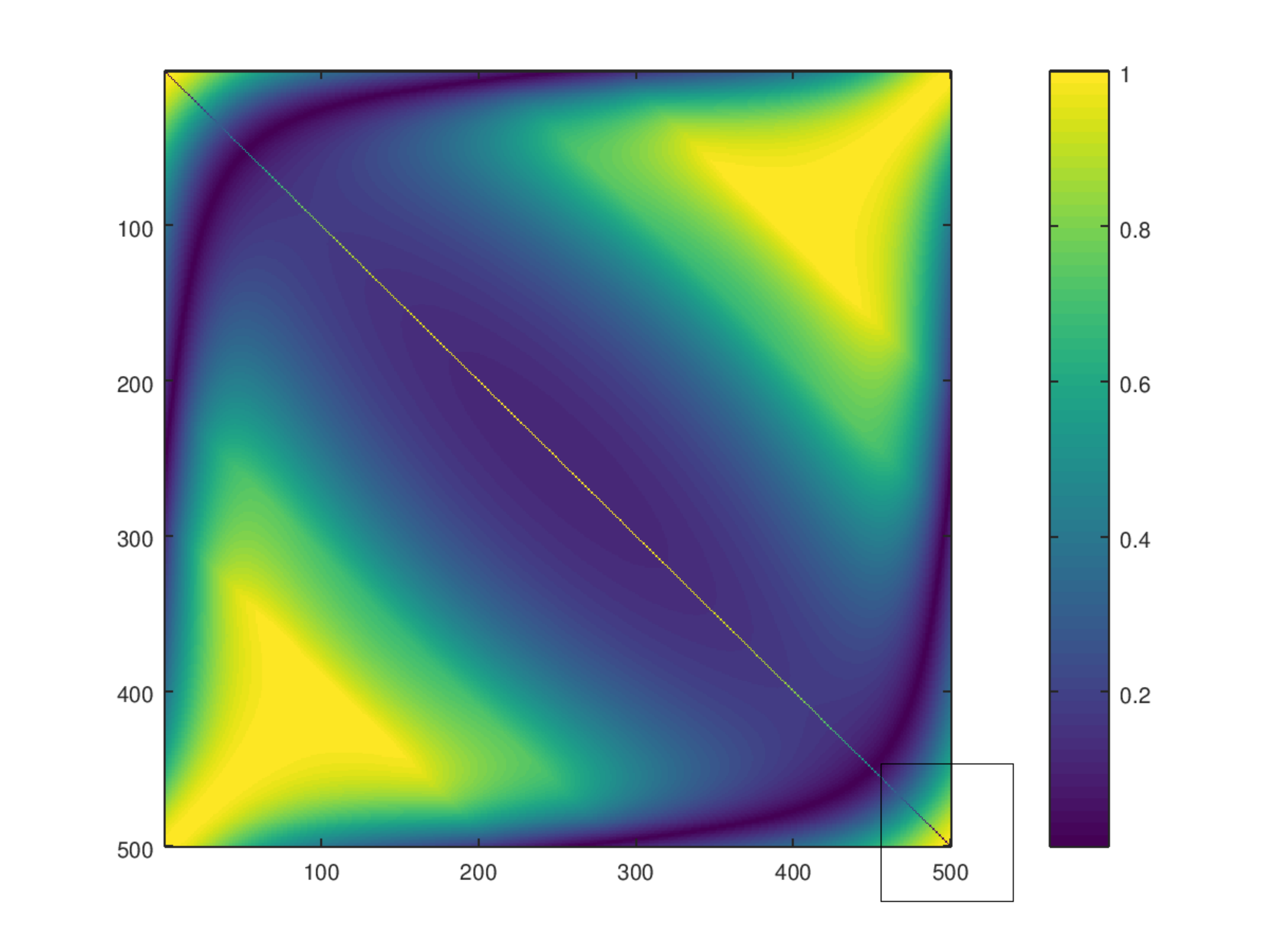}
		}
	\parbox{7cm}{
		\includegraphics[width = 7 cm]{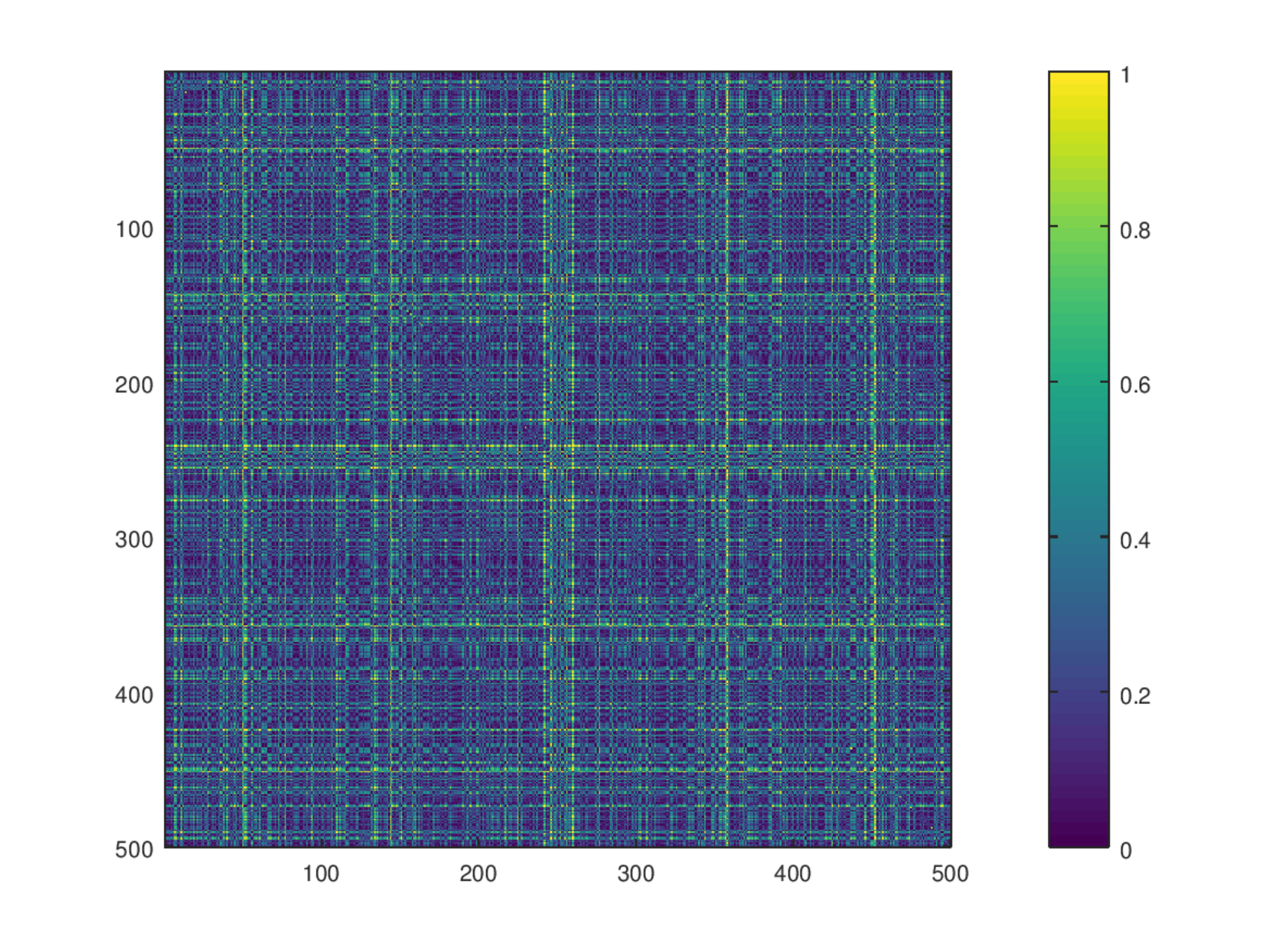}
		}
		\caption{ Cycle-wise dominant entries of the matrix $B = WAW^{\dagger}$. Heat map represents the magnitude of the matrix entry normalized by the largest magnitude among entries in the corresponding cycle. The matrices $A$ are symmetric Toeplitz matrices of dimension 500. Left: When $A$ has  entries from example-1 of \cite{cai2005generalization}, and note the higher relative magnitudes of the cycles at the bottom right indicated by a box. Right: When $A$ has random Gaussian $\mathcal{N}(0,1)$ entries. }\label{heatmap_figure}
\end{figure}

Toeplitz linear systems appear, for example, in solving linear ordinary differential equations and delay differential equations. Here, T. Chan and generalized T. Chan's preconditioners are used in speeding up Conjugate Gradient (PCG) algorithms \cite{cai2005generalization, jin2008survey}. In example-1 of the article \cite{cai2005generalization}, preconditioned CG is applied on a symmetric Toeplitz matrix with first row given by
    $A(1,:) = \begin{bmatrix} 
    2 &\frac{-1}{2} & \frac{-1}{2^2} & \frac{-1}{2^3} & \cdots & \frac{-1}{2^{n-1}}
    \end{bmatrix}$
    with right hand side vector $b = \begin{bmatrix}
    1 \\
    2\\
    3 \\
    \vdots \\ 
    n
    \end{bmatrix}$. The preconditioner considered was a genaralization of the T. Chan's preconditioner, given as
    $P(k) = W^{\dagger} (B \circ Q(k)) W$ with the matrix $Q(k)$ given as
    $Q(k) = \begin{bmatrix}
    I_{n-\ceil{\sqrt{n-k}}} & 0 \\
    0 & \mathbf{1}_{\ceil{\sqrt{n-k}}}
    \end{bmatrix}$
    Where $\mathbf{1}_k$ is a $k \times k$ matrix of all ones.
    Note that this preconditioner includes a single circulant component i.e. the diagonal of $WAW^{\dagger}$ along with a square sub-matrix at the bottom-right corner in the matrix $B$. Table \ref{t_it} shows the number of iterations required by a preconditioned CG for positive definite symmetric Toeplitz matrices. The generalized T. Chan's preconditioner performs marginally better than the proposed banded preconditioner with a few cycles of $WAW^{\dagger}$, for example-1, and not in the other cases. This can be explained using the Figure \ref{heatmap_figure}. The matrix $WAW^{\dagger}$ for example-1 has higher relative magnitudes of all cycles in the bottom right corner. As mentioned, the T. Chan's preconditioner can capture most of these entries as a bottom-right sub-matrix, in comparison with the corresponding banded preconditioner. Whereas other Toeplitz matrices may not show such a cluster of high magnitude entries near the diagonal. Thus, T. Chan's preconditioner does not reduce required iterations with small increments in the size of its non-zero sub-matrix, in general. Thus the banded preconditioner using the cycles of $WAW^{\dagger}$ performs better than the genalized T. Chan's preconditioner in the other positive definite Toeplitz matrices with arbitrary entries. 
    
    In the case of block-Toeplitz matrices we see that the banded preconditioner speeds up computations with inclusion of cycles of $WAW^{\dagger}$ as shown in Table \ref{bt_it}, whereas the T. Chan's preconditioner fails to do so for the same number of non-zero entries in it. These results highlight the generality of the dominant cycles of matrix $B=WAW^{\dagger}$, as an appropriate preconditioner for matrices with some periodicity in entries.
    
\begin{table}
    \begin{tabular}{|m{4cm}|m{3em}|m{3em}|m{3em}|m{3em}|m{3em}|m{3em}|} \hline
    Number of non-zero entries in preconditioner and iterations required $\longrightarrow$& I & P(n) & P(3n) & P(5n) & P(7n) & P(9n) \\ \hline
    Generalized T. Chan &   683& 30& 23& 23& 23& 23 \\ \hline
    $WAW^{\dagger}$ cycles &  &  30&44&43&45&47 \\ \hline
     Generalized T. Chan&  61&38&38&38&38&38 \\ \hline
    $WAW^{\dagger}$ cycles &  &  38&27&24&21&21\\ \hline
     Generalized T. Chan &  103& 61& 61& 61& 61& 61\\ \hline
    $WAW^{\dagger}$ cycles &  & 61&45&40&34&33 \\ \hline
     Generalized T. Chan&  51&33&33&33&33&33\\ \hline
    $WAW^{\dagger}$ cycles &  &   33&23&21&21&19\\ \hline
     Generalized T. Chan &  56&36&36&36&36&34 \\ \hline
    $WAW^{\dagger}$ cycles &  & 36&26&23&22&20\\ \hline
    \end{tabular}
    \caption{Toeplitz matrices of dimension 2000. The first case is example-1 of \cite{cai2005generalization}, and others were generated using random entries. The preconditioner `I' represents the plain solver. Both preconditioners are identical when only a single cycle of $WAW^{\dagger}$ is considered. The stopping criterion was a relative residue less than $10^{-6}$.}\label{t_it}
   \end{table}

   \begin{table}
    \begin{tabular}{|m{4cm}|m{3em}|m{3em}|m{3em}|m{3em}|m{3em}|m{3em}|m{3em}|} \hline
    Number of non-zero entries in preconditioner and iterations required $\longrightarrow$& I & P(n) & P(3n) & P(5n) & P(7n) & P(9n) & P(11n)\\ \hline
    Generalized T. Chan& 168&162&162&162&162&162&162  \\ \hline
    $WAW^{\dagger}$ cycles &  & 162&148&121&106& 78& 19 \\ \hline
    Generalized T. Chan&136&121&121&121&121&121&121 \\ \hline
    $WAW^{\dagger}$ cycles &  &121&119&106& 94& 66& 17\\ \hline
    Generalized T. Chan& 172&172&172&172&172&172&172\\ \hline
    $WAW^{\dagger}$ cycles & &172&131&134&105& 66& 19\\ \hline
    Generalized T. Chan&408&402&403&402&400&403&399 \\ \hline
    $WAW^{\dagger}$ cycles & & 402&375&327&284&232& 29 \\ \hline
    Generalized T. Chan&125&124&124&124&124&124&124 \\ \hline
    $WAW^{\dagger}$ cycles &  & 124&120&100& 90& 51& 17\\ \hline
    \end{tabular}
     \caption{Block-Toeplitz matrices of dimension 1100 of block size 11. The matrices are symmetrix positive definite with symmetric blocks, and condition number of the matrices is $\mathcal{O}(10^4)$. The stopping criterion was a relative residue less than $10^{-6}$. The preconditioner `I' represents the plain solver. Both preconditioners are identical when only a single cycle of $WAW^{\dagger}$ is considered.}\label{bt_it}
    \end{table}

\section{Conclusion}
We began with a decomposition of any given matrix into circulant matrices with periodic relaxations on the unit circle. Exploiting the periodicity of entries along the diagonals, and the dominance of a few circulant components, we can reduce the given matrix using appropriate fast-Fourier-transform operations for an approximate and sparse similarity transformation. In the appendix, we also compare the efficacy of the symbol of a Toeplitz matrix, and its one-term circulant approximation, in describing the spectra. While both require $\mathcal{O}(n^2)$ operations each for the evaluation, the latter is more generally applicable to matrices with any periodicity in entries.

Using numerical results, we highlighted the efficacy of the sparsification of the periodic matrices using the circulant decomposition, in terms of errors both in approximation of entries of the matrix and also its eigenvalues. Examples of the relative errors in the eigenvalues were produced as a function of the number of circulant components included in the approximation. Results in preconditioning linear systems were presented where the generalized T. Chan's preconditioner was compared. These results show that the suggested sparse similarity transformation of a matrix is useful in efficiently approximating eigenvalues, and preconditioning linear systems, and may as well be exploited for other evaluations when a dense matrix has any periodicity along its diagonals.

\appendix

\section{Single-term circulant approximation of a Toeplitz matrix, and its symbol}
A Toeplitz matrix is of the form,
\vspace{3mm}
\newline
$A_n = \begin{bmatrix}
a_0 & a_1 & a_2 &  \cdots & a_{n-1} \\
a_{-1} & a_0 & a_1 & a_2  & \ddots \\
a_{-2}& a_{-1} & a_0 & \ddots  & \ddots \\
\vdots & \ddots & \ddots    & \ddots & \ddots \\
a_{-(n-1)} & a_{-(n-2)}& \cdots & \cdots &a_0 \\
\end{bmatrix}$. \\

The function $a(e^{i \theta}) = \sum \limits_{k = -\infty}^{\infty}a_k e^{ i k \theta}$ for $ \theta \in [0,2 \pi)$ is called the \emph{symbol} of the family of Toeplitz matrices $A_n$ \cite{bottcher2012introduction}. This section considers different types of symbols, the corresponding eigenvalues of the single-term circulant approximation, and the actual spectrum of the Toeplitz matrices. The single-term approximation is given by only the diagonal entries of $WA_nW^{\dagger}$, and is closely related to the symbol of the Toeplitz matrices. Both, the former and the latter, allow approximation of eigenvalues of Toeplitz matrices in $\mathcal{O}(n^2)$ arithmetic operations. But the circulant approximation with one or more terms is more generally applicable even to block-Toeplitz and other dense matrices with some periodicity in entries.  

\subsection*{Case 1: the symbol is a product of a polynomial and a trigonometric polynomial} 
\begin{figure}
	\begin{center}
		\includegraphics[width = 9 cm]{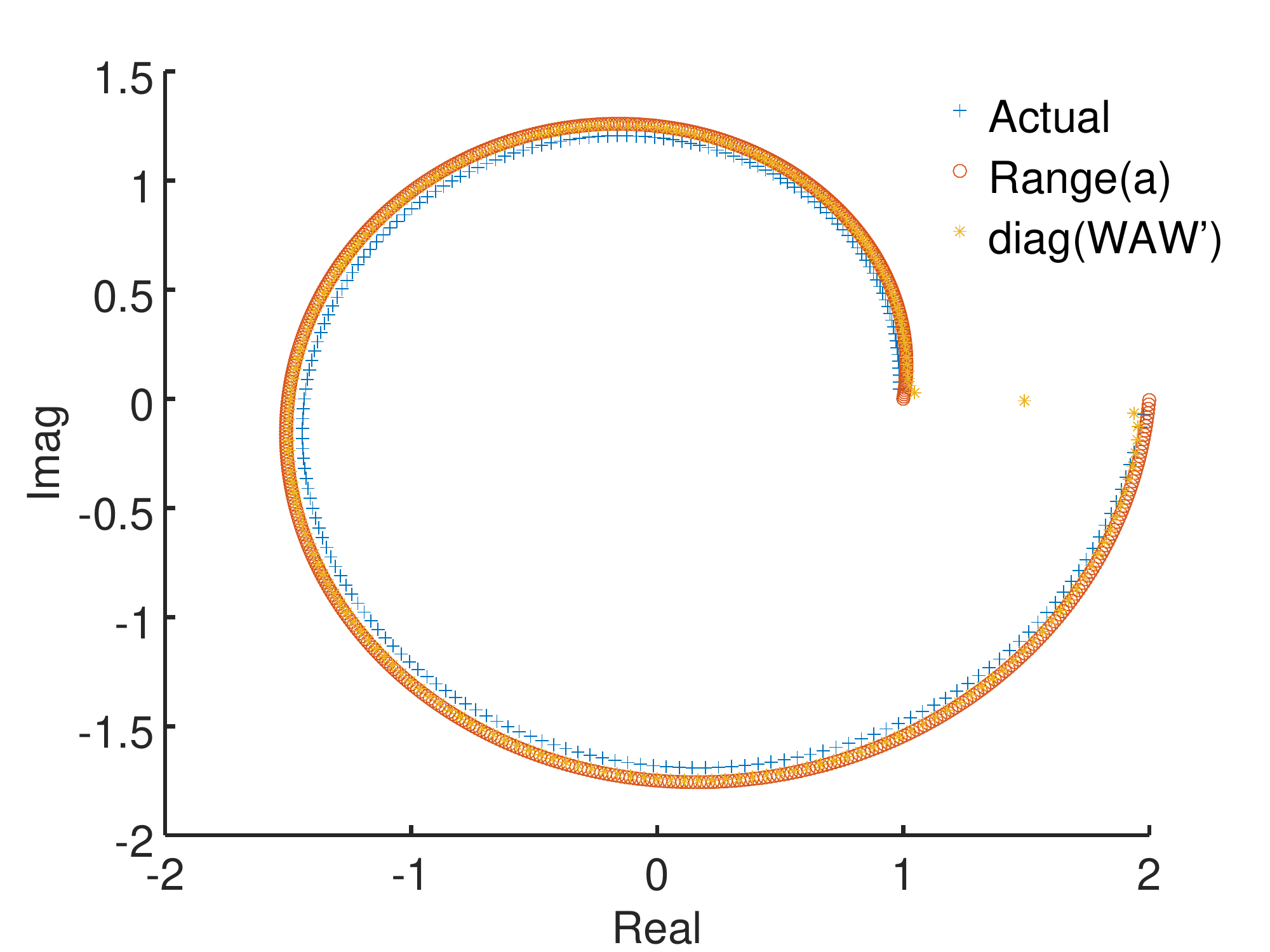}
		\caption{Range of the symbol, and diagonal entries of $WA_nW^{\dagger}$ for a Toeplitz matrix of dimension 200 with a symbol $a(\theta) = (1+ \theta) e^{i \theta}$.}\label{symbolpq}
	\end{center}
\end{figure}

Spectrum of $A_n$ is said to be canonically distributed if the limiting spectrum approaches the range of the symbol. For symbols of the form $a(e^{i \theta}) = p(\theta)q(e^{i \theta})$ with polynomials $p$ and $q$, we consider the Toeplitz matrices whose entries are from its truncated Fourier series expansion. Several conditional theorems on the symbol such that the spectrum of $A_n$ is canonically distributed are presented in \cite{bottcher2012introduction}. It was shown that if the complement of the range is a connected set, then the spectrum follows the symbol \cite{tilli1999some}. Figure \ref{symbolpq}  shows the spectrum, single-term circulant approximation, and the symbol for $a( \theta) = (1+ \theta) e^{i \theta}$.

\subsection*{Case 2: the symbol of the banded Toeplitz matrix is a trigonometric polynomial} 
\begin{figure}
	\begin{center}
		\includegraphics[width = 9 cm]{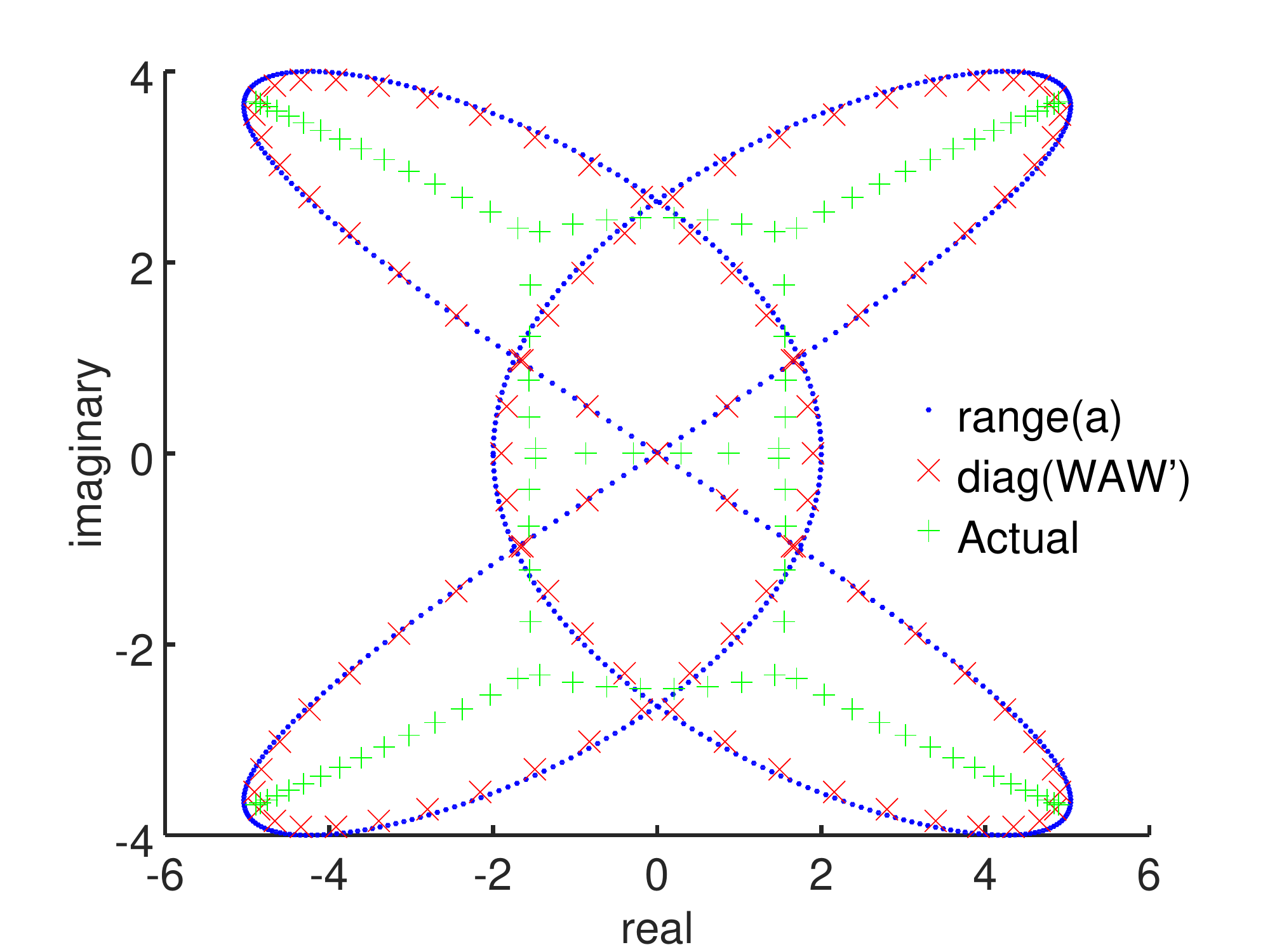}
		\caption{Range of the symbol, and diagonal entries of $WA_nW^{\dagger}$ for a banded Toeplitz matrix of dimension 100, with 7 non-zero entries in each row. The work in this paper is not directed at such matrices.}\label{symbol1}
	\end{center}
\end{figure}

For a banded Toeplitz matrix $A_n$, with first row $a_0,a_1,a_2 , \cdots a_l$ and first column entries $a_0,a_{-1}, \cdots a_{-m}$, the symbol given by $a(e^{i \theta}) = \sum \limits_{k = -m}^{l}a_k e^{i k \theta}$ is a curve in the complex plane. The eigenvalues of such matrices lie in the convex hull of the curve $a(e^{i \theta})$ for $\theta \in [0,2 \pi)$ \cite{bottcher2005spectral, tilli1999some}. On the other hand, we have the diagonal entries of $WA_nW^{\dagger}$ as FFT of the sequence $$(a_0, \frac{n-1}{n}a_1, \frac{n-2}{n}a_2, \cdots \frac{n-l}{n}a_l, 0 ,0 , \cdots 0 , \frac{n-m}{n}a_{-m}, \cdots \frac{n-1}{n}a_{-1}).$$ Note that these values lie in near proximity to the trace of the symbol in the complex plane when $k \ll n$, as shown in Figure \ref{symbol1}.

\subsection*{Case 3: the symbol is a sum of a polynomial and a trigonometric polynomial} 
\begin{figure}
	\begin{center}
	   \includegraphics[width = 9 cm]{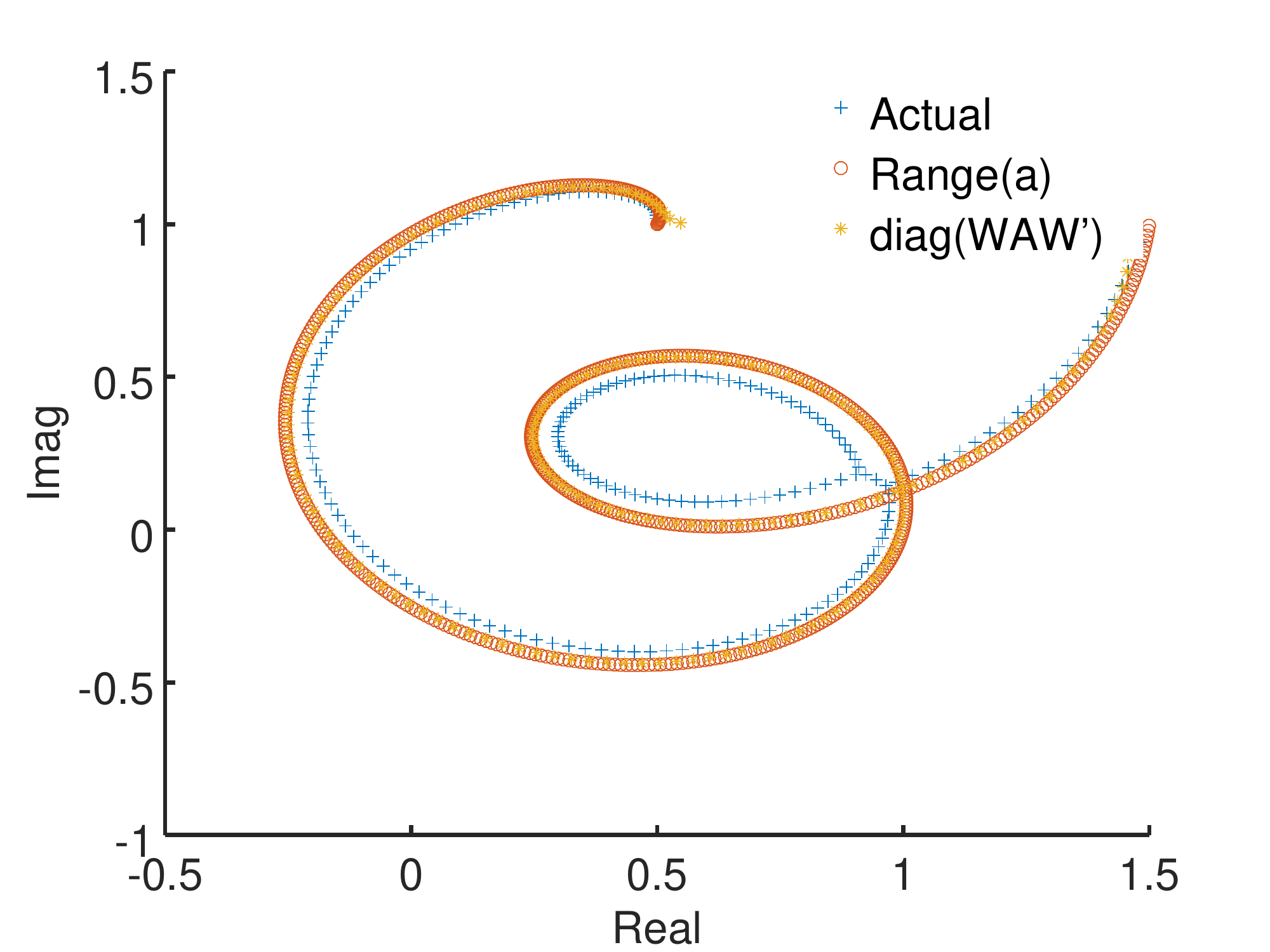}
		\caption{Range of the symbol, and diagonal entries of $WA_nW^{\dagger}$ for a Toeplitz matrix of dimension 200 with a symbol $a(\theta) = \frac{\theta}{2 \pi} + i \frac{1}{\pi^2}(\theta - \pi)^2 + e^{i 2 \theta}$.}\label{symbolp+ip+q}
	\end{center}
\end{figure}

The spectrum of the matrix with symbol $a(e^{i \theta}) = \frac{\theta}{2 \pi} + i \frac{1}{\pi^2}(\theta - \pi)^2 + e^{i 2 \theta} $ is shown in Figure \ref{symbolp+ip+q}. This can be regarded as the combination of the previous two cases, even though the spectra mimic the symbol, they lie inside the convex hull.

\FloatBarrier
\noindent \textbf{Author declarations}:\\

\noindent \textbf{Funding}: Murugesan Venkatapathi acknowledges the support of the Science and Engineering Research Board (SERB) grant MTR/2019/000712 in performing this research.\\
\vspace{2mm}

\noindent \textbf{Conflicts of interest}: The authors do not have any competing financial or non-financial interests to declare.

\bibliographystyle{unsrt}
\bibliography{ref.bib}

\begin{thebibliography}{10}

\bibitem{Ambikasaran2013Dense}
Sivaram Ambikasaran and Eric Darve.
\newblock An {$O(N$log$N)$} fast direct solver for partial hierarchically
  semi-separable matrices.
\newblock {\em Journal of Scientific Computing}, 57(3):477--501, 2013.

\bibitem{saad2011numerical}
Yousef Saad.
\newblock {\em Numerical methods for large eigenvalue problems: revised
  edition}, volume~66.
\newblock SIAM, 2011.

\bibitem{koslowski1993linear}
Th~Koslowski and W~Von~Niessen.
\newblock Linear combination of {L}anczos vectors: A storage-efficient
  algorithm for sparse matrix eigenvector computations.
\newblock {\em Journal of Computational Chemistry}, 14(7):769--774, 1993.

\bibitem{paige1971computation}
Christopher~Conway Paige.
\newblock {\em The computation of eigenvalues and eigenvectors of very large
  sparse matrices.}
\newblock PhD thesis, University of London, 1971.

\bibitem{bunch2014sparse}
James~R Bunch and Donald~J Rose.
\newblock {\em Sparse matrix computations}.
\newblock Academic Press, 2014.

\bibitem{ekstrom2018sparse}
Sven-Erik Ekström, Carlo Garoni, and Stefano Serra-Capizzano.
\newblock Are the eigenvalues of banded symmetric {T}oeplitz matrices known in
  almost closed form?
\newblock {\em Experimental Mathematics}, 27(4):478--487, 2018.

\bibitem{hariprasad2021chain}
M~Hariprasad and Murugesan Venkatapathi.
\newblock Semi-analytical solutions for eigenvalue problems of chains and
  periodic graphs.
\newblock {\em Applied Mathematics and Computation}, 411:126512, 2021.

\bibitem{trench1989numerical}
F~William Trench.
\newblock Numerical solution of the eigenvalue problem for {H}ermitian
  {T}oeplitz matrices.
\newblock {\em SIAM Journal on Matrix Analysis and Applications},
  10(2):135--146, 1989.

\bibitem{trench1985eigenvalue}
William~F Trench.
\newblock On the eigenvalue problem for {T}oeplitz band matrices.
\newblock {\em Linear Algebra and its Applications}, 64:199--214, 1985.

\bibitem{bini1988efficient}
D~Bini and V~Pan.
\newblock Efficient algorithms for the evaluation of the eigenvalues of (block)
  banded {T}oeplitz matrices.
\newblock {\em Mathematics of Computation}, 50(182):431--448, 1988.

\bibitem{luk2000fast}
Franklin~T Luk and Sanzheng Qiao.
\newblock A fast eigenvalue algorithm for {H}ankel matrices.
\newblock {\em Linear Algebra and its Applications}, 316(1-3):171--182, 2000.

\bibitem{Gray1972Toeplitz}
R.~Gray.
\newblock On the asymptotic eigenvalue distribution of {T}oeplitz matrices.
\newblock {\em IEEE Transactions on Information Theory}, 18(6):725--730, 1972.

\bibitem{Boo1998Toeplitz}
N.K. Bose and K.J. Boo.
\newblock Asymptotic eigenvalue distribution of block-{T}oeplitz matrices.
\newblock {\em IEEE Transactions on Information Theory}, 44(2):858--861, 1998.

\bibitem{bottcher2005spectral}
Albrecht B{\"o}ttcher and Sergei~M Grudsky.
\newblock {\em Spectral properties of banded {T}oeplitz matrices}.
\newblock SIAM, 2005.

\bibitem{bottcher2012introduction}
Albrecht B{\"o}ttcher and Bernd Silbermann.
\newblock {\em Introduction to large truncated {T}oeplitz matrices}.
\newblock Springer Science \& Business Media, 2012.

\bibitem{zhu2017asymptotic}
Zhihui Zhu and Michael~B Wakin.
\newblock On the asymptotic equivalence of circulant and {T}oeplitz matrices.
\newblock {\em IEEE Transactions on Information Theory}, 63(5):2975--2992,
  2017.

\bibitem{bogoya2022fast}
Manuel Bogoya, Sven-Erik Ekstr{\"o}m, and Stefano Serra-Capizzano.
\newblock Fast {T}oeplitz eigenvalue computations, joining
  interpolation-extrapolation matrix-less algorithms and simple-loop theory.
\newblock {\em Numerical Algorithms}, pages 1--24, 2022.

\bibitem{chan1987asymptotic}
Raymond Chan and Gilbert Strang.
\newblock The asymptotic {T}oeplitz-circulant eigenvalue problem.
\newblock 1987.

\bibitem{arbenz1991computing}
Peter Arbenz.
\newblock Computing eigenvalues of banded symmetric {T}oeplitz matrices.
\newblock {\em SIAM Journal on Scientific and Statistical Computing},
  12(4):743--754, 1991.

\bibitem{chan1988optimal}
Tony~F Chan.
\newblock An optimal circulant preconditioner for {T}oeplitz systems.
\newblock {\em SIAM Journal on Scientific and Statistical Computing},
  9(4):766--771, 1988.

\bibitem{jin2008survey}
Xiao-Qing Jin and Yi-Min Wei.
\newblock A survey and some extensions of {T. Chan’s} preconditioner.
\newblock {\em Linear Algebra and its Applications}, 428(2-3):403--412, 2008.

\bibitem{Reichel2012circulant}
Silvia Noschese and Lothar Reichel.
\newblock Generalized circulant strang-type preconditioners.
\newblock {\em Numerical Linear Algebra with Applications}, 19(1):3--17, 2012.

\bibitem{Noschese2014circulant}
Silvia Noschese and Lothar Reichel.
\newblock A note on superoptimal generalized circulant preconditioners.
\newblock {\em Applied Numerical Mathematics}, 75:188--195, 2014.

\bibitem{cai2005generalization}
Ming-Chao Cai, Xiao-Qing Jin, and Yi-Min Wei.
\newblock A generalization of {T. Chan’s} preconditioner.
\newblock {\em Linear algebra and its applications}, 407:11--18, 2005.

\bibitem{achlioptas2007fast}
Dimitris Achlioptas and Frank McSherry.
\newblock Fast computation of low-rank matrix approximations.
\newblock {\em Journal of the ACM (JACM)}, 54(2):9--es, 2007.

\bibitem{tilli1999some}
Paolo Tilli.
\newblock Some results on complex {T}oeplitz eigenvalues.
\newblock {\em Journal of Mathematical Analysis and Applications},
  239(2):390--401, 1999.

\end{thebibliography}
\end{document}